\documentclass{aptpub}
\usepackage{hyperref}
\usepackage[numbers,square]{natbib}

\newcommand{\bS}{\mathbb{S}}

\newcommand{\cF}{\mathcal{F}}

\newcommand{\cQ}{\mathcal{Q}}
\newcommand{\cP}{\mathcal{P}}

\DeclareMathOperator{\relint}{relint}

\newcommand{\boldE}{\mathbb E}

\newcommand{\R}{\mathbb{R}}
\newcommand{\boldN}{\mathbb{N}}

\newcommand{\C}{\mathbb{C}}

\renewcommand{\P}{\mathbb{P}}

\renewcommand{\Re}{\operatorname{Re}}

\newcommand{\erf}{\operatorname{erf}}

\newcommand{\simn}{\underset{n \to \infty}{\sim}}

\newcommand{\Var}{\mathop{\mathrm{Var}}\nolimits}
\newcommand{\Vol}{\mathop{\mathrm{Vol}}\nolimits}

\newcommand{\Cov}{\mathop{\mathrm{Cov}}\nolimits}

\newcommand{\conv}{\mathop{\mathrm{Conv}}\nolimits}
\newcommand{\pos}{\mathop{\mathrm{pos}}\nolimits}

\newcommand{\eps}{\varepsilon}
\newcommand{\eqdistr}{\stackrel{d}{=}}

\newcommand{\bsl}{\backslash}

\newcommand{\dd}{{\rm d}}
\newcommand{\eee}{{\rm e}}

\authornames{Zakhar Kabluchko and Dmitry Zaporozhets} 
\shorttitle{Gaussian polytopes and regular spherical simplices} 

\numberwithin{equation}{section}  

\begin{document}

\title{Absorption probabilities for Gaussian polytopes\\ and regular spherical simplices}

\authorone[Westf\"alische Wilhelms-Universit\"at M\"unster]{Zakhar Kabluchko}
\addressone{
Orl\'eans--Ring 10,
48149 M\"unster, Germany}

\authortwo[St.\ Petersburg Department of Steklov Mathematical Institute]{Dmitry Zaporozhets}
\addresstwo{Fontanka~27, 191011 St.\ Petersburg, Russia}

\begin{abstract}
The Gaussian polytope $\mathcal P_{n,d}$ is the convex hull of $n$ independent standard normally distributed points in $\R^d$. We derive explicit expressions for the probability that $\mathcal P_{n,d}$ contains a fixed point $x\in\R^d$ as a function of the Euclidean norm of $x$, and the probability that $\mathcal P_{n,d}$ contains the point $\sigma X$, where $\sigma\geq 0$ is constant and $X$ is a standard normal vector independent of $\mathcal P_{n,d}$. As a by-product, we also compute the expected number of $k$-faces and the expected volume of $\mathcal P_{n,d}$, thus recovering the results of Affentranger and Schneider [\textit{Discr.\ and Comput. Geometry}, 1992] and  Efron [\textit{Biometrika}, 1965], respectively. All formulas are in terms of the volumes of regular spherical simplices, which, in turn, can be expressed through the standard normal distribution function $\Phi(z)$ and its complex version $\Phi(iz)$. The main tool used in the proofs is the conic version of the Crofton formula.
\end{abstract}



\section{Statement of main results}\label{1306}

\subsection{Introduction}
Let $X_1,\ldots,X_n$ be independent random vectors with standard Gaussian distribution on $\R^d$.
The \textit{Gaussian polytope} $\cP_{n,d}$ is defined as the convex hull of $X_1,\ldots,X_n$, that is
$$
\cP_{n,d}
=
\conv(X_1,\ldots,X_n) =
\left\{\sum_{i=1}^n \lambda_i X_i \colon \lambda_1,\ldots,\lambda_n \geq 0,  \sum_{i=1}^n \lambda_i = 1\right\}.
$$
The main aim of the present paper is to provide an explicit expression for the \textit{absorption probability}, that is the probability that $\cP_{n,d}$ contains a given deterministic point $x\in\R^d$.
By rotational symmetry,  the absorption probability depends only on the Euclidean norm $|x|$. It turns out that it is more convenient to pass to the complementary event and consider the \textit{non-absorption probability}
\begin{equation}\label{eq:def_f}
f_{n,d}(|x|) := \P[x\notin \conv (X_1,\ldots,X_n)].
\end{equation}
A classical result of Wendel~\cite{wendel} (which is valid in a setting more general than the Gaussian one considered here), see also~\cite[Theorem 8.2.1]{SW08},   states that
\begin{equation}\label{eq:wendel}
f_{n,d}(0) = \frac1{2^{n-1}} \left(\binom{n-1}{d-1} + \binom{n-1}{d-2} +\ldots \right).
\end{equation}
We shall compute $f_{n,d}(r)$ for general $r\geq 0$. The main idea is that we shall make the point $x$ random, with a rotationally symmetric Gaussian distribution and certain variance $\sigma^2\geq 0$.
Namely, let $X$ be a $d$-dimensional standard Gaussian random vector which is independent of $X_1,\ldots,X_n$. We shall compute
\begin{equation}\label{eq:eq:p_n_d_def}
p_{n,d} (\sigma^2) := \P[\sigma X \notin \conv (X_1,\ldots, X_n)].
\end{equation}
This probability can be related to a certain Laplace-type transform of $f_{n,d}$. After inverting the Laplace transform, we shall obtain a formula for
$f_{n,d}$. This formula involves certain function $g_n(r)$ which expresses the volume of regular spherical simplices and which will be studied in detail below.

The probability $p_{n,d}(\sigma^2)$ is closely related to the expected number of faces of the polytope $\cP_{n,d}$. Let $f_k(\cP_{n,d})$ be the number of $k$-dimensional faces ($k$-faces) of $\cP_{n,d}$. Exact formulas for $\boldE f_{k}(\cP_{n,d})$ were derived by~\citet[\S 4]{renyi_sulanke1} (for $d=2$), \citet{efron} (for $d=2,3$), \citet{raynaud} (for faces of maximal dimension, that is for  $k=d-1$). \citet{affentranger} proved an asymptotic formula valid for general $d$ and $k$; see also \citet{carnal} for the case $d=2$. \citet{baryshnikov_vitale} showed that the expected number of $k$-faces of $\cP_{n,d}$ is the same as the expected number of $k$-faces of a random projection of the regular simplex with $n$ vertices onto a uniformly chosen linear subspace of dimension $d$ (the so-called \textit{Goodman--Pollack model}). Finally, \citet{AS92} expressed the expected number of $k$-faces of the random projection of any polytope in terms of the internal and external angles of that polytope. Combining the results of \citet{AS92} and \citet{baryshnikov_vitale}, one obtains an expression for $\boldE f_k(\cP_{n,d})$ in terms of the internal and external angles of the regular simplex.

\citet{HMR04} expressed some important functionals of the Gaussian polytope including the expected number of $k$-faces through the probabilities of the form $p_{n,d}(\sigma^2)$ and computed their asymptotics. As a by-product of our results, we shall provide explicit formulas for these functionals, thus recovering the results obtained in~\cite{AS92} and~\cite{HMR04}.
Recent surveys on random polytopes can be found in~\cite{hug_survey,MCR10,schneider_polytopes}.

\subsection{Non-absorption probabilities}\label{subsec:main_results}
Our explicit formulas will be stated in terms of the functions $g_n(r)$, $n\in\boldN_0$, defined by $g_0(r):=1$ and
\begin{equation}\label{eq:def_g}
g_n(r) = \P[\eta_1<0,\ldots,\eta_n<0], \quad r\geq -1/n, \;\; n\in\boldN,
\end{equation}
where  $(\eta_1,\ldots,\eta_n)$ is a zero-mean Gaussian vector with
\begin{equation}\label{eq:cov_eta}
\Cov (\eta_i,\eta_j)
=
\begin{cases}
1+r, &\text{ if } i=j,\\
r, &\text{ if } i\neq j.
\end{cases}
\end{equation}
The fact that~\eqref{eq:cov_eta} indeed defines a valid (i.e.,\ positive semi-definite) covariance matrix for $r\geq -1/n$ is easily verified using the inequality between the arithmetic and quadratic means.

Many known and some new properties of the function $g_n$ (which is closely related to the Schl\"afli function~\cite{boehm_hertel_book}) will be collected in Sections~\ref{sec:geometric_g_n} and~\ref{sec:explicit_g_n}.
At this place, we just state an explicit expression for $g_n$ in terms of the standard normal distribution function $\Phi$.
It is known that $\Phi$ admits an analytic continuation to the entire complex plane. We shall need its values on the real and imaginary axes, namely
\begin{equation}\label{eq:Phi_def}
\Phi(z) = \frac{1}{\sqrt{2\pi}} \int_{-\infty}^z \eee^{-t^2/2} \dd t,
\quad
\Phi( i z) = \frac 12 + \frac i{\sqrt{2\pi}} \int_0^{z} \eee^{t^2/2} \dd t,
\quad z\in\R.
\end{equation}
The reader more used to the error function $\erf$ may transform everything by applying  the formula $\Phi(z) = 1/2 + \erf(z/\sqrt 2)$. With this notation, an explicit formula for $g_n$ reads as follows:
\begin{align}
g_n(r)
&=
\frac 1 {\sqrt {2\pi}} \int_{-\infty}^{\infty} \Phi^n (\sqrt r x) \eee^{-x^2/2} \dd x \label{eq:g_n_first_formula1}\\
&=
\begin{cases}
\frac {1}{\sqrt{2\pi}} \int_{-\infty}^{\infty} \left(\frac{1}{\sqrt{2\pi}} \int_{-\infty}^{x\sqrt r} \eee^{-z^2/2} \dd z\right)^{n}\eee^{-x^2/2}\dd x, &\text{if } r\geq 0,\\
\frac {2}{\sqrt{\pi}} \int_0^{\infty} \Re \left[\left(\frac 12 + \frac{i}{\sqrt{\pi}} \int_0^{x\sqrt{-r}} \eee^{z^2} \dd z\right)^{n}\right]\eee^{-x^2}\dd x, &\text{if } -1/n\leq r \leq 0,\label{eq:g_n_first_formula2}
\end{cases}
\end{align}
where in~\eqref{eq:g_n_first_formula1} we agree that $\sqrt{r} = i\sqrt{-r}$ if $r<0$. The next theorem provides a formula for the probability that $\sigma X \notin \cP_{n,d}$.
\begin{thm}\label{theo:main1}
Let $X,X_1,\ldots,X_n$ be independent standard Gaussian random vectors in $\R^d$, where $n\geq d+1$.  Then, for every $\sigma\geq 0$,
\begin{equation}
p_{n,d} (\sigma^2) = \P[\sigma X \notin \conv (X_1,\ldots, X_n)] = 2(b_{n,d-1}(\sigma^2) + b_{n,d-3}(\sigma^2) +\ldots),
\end{equation}
where
\begin{equation}\label{eq:def_b_k}
b_{n,k}(\sigma^2) = \binom nk g_k\left(-\frac{\sigma^2}{1+k\sigma^2}\right) g_{n-k} \left(\frac{\sigma^2}{1+k\sigma^2}\right),
\end{equation}
for $k=0,\ldots,n$, and $b_{n,k}=0$ for $k\notin \{0,1,\ldots,n\}$.
\end{thm}
The proof of Theorem~\ref{theo:main1} will be given in Section~\ref{sec:proof_main}. The main idea is to interpret $p_{n,d}(\sigma^2)$ as the probability that a uniform  random linear subspace intersects certain $n$-dimensional convex cone $C= C_n(\sigma^2)$.
By the conic Crofton formula (which will be recalled in Theorem~\ref{theo:crofton_conic} below), this intersection probability can be expressed in terms of the conic intrinsic volumes $\upsilon_0(C),\ldots, \upsilon_n(C)$ of $C$.
At this point, we can forget about the original problem and concentrate on computing the conic intrinsic volumes, which is a purely geometric problem. It turns out that $\upsilon_k(C) = b_{n,k}(\sigma^2)$; see Proposition~\ref{prop:upsilon_n_k}, below.

\begin{ex}[Wendel's formula]
Let $\sigma^2 = 0$. By symmetry reasons it is clear (and will be stated in Proposition~\ref{prop:g_n_properties} (d)) that $g_k(0) = 2^{-k}$, $g_{n-k}(0) = 2^{-(n-k)}$ and hence $b_{n,k}(0) = 2^{-n} \binom nk$. Theorem~\ref{theo:main1} simplifies to
\begin{align}
\P[0\notin \conv (X_1,\ldots,X_n)]
&=
\frac1{2^{n-1}} \left(\binom{n}{d-1} + \binom{n}{d-3} +\ldots \right) \label{eq:wendel_formula_rep}\\
&=
\frac1{2^{n-1}} \left(\binom{n-1}{d-1} + \binom{n-1}{d-2} +\ldots \right) \notag
,
\end{align}
where in the second line we used the defining property of the Pascal triangle.
This recovers Wendel's formula~\eqref{eq:wendel} in the Gaussian case.
\end{ex}

By conditioning on $X$ in Theorem~\ref{theo:main1}, we shall derive the following
\begin{cor}\label{cor:intensity}
The function $f_{n,d}(|x|) = \P[x \notin \conv(X_1,\ldots,X_n)]$ satisfies
$$
\int_0^\infty f_{n,d}(\sqrt{2u}) u^{\frac d2 - 1} \eee^{-\lambda u} \dd u = 2 \Gamma (d/2) \lambda^{-\frac d2}  (b_{n,d-1} (1/\lambda) + b_{n,d-3} (1/\lambda) + \ldots)
$$
for all $\lambda>0$.
\end{cor}

It is possible to  invert the Laplace transform explicitly. Recall that $\Phi$ is the standard normal distribution function.

\begin{thm}\label{theo:main_intensity}
For all $u> 0$ we have
$$
f_{n,d}(\sqrt{2u}) - f_{n,d}(0) =  2 u^{1-(d/2)} (a_{n,d-1}(u) + a_{n,d-3}(u) + \ldots),
$$
where $f_{n,d}(0)$ is given by Wendel's formula~\eqref{eq:wendel} and
\begin{align}
a_{n,k}(u) &= \binom nk \int_0^u \eee^{-vk} F'_{k,n-k} (v) (u-v)^{(d/2)-1} \dd v,\label{eq:a_k} \\
F_{k,n-k} (v) &=\frac 1 {\pi}
\int_0^v
\left(\frac{\Phi^{n-k}(\sqrt{2 w}) + \Phi^{n-k}(-\sqrt{2 w})}{2\sqrt w} \right.
\notag
\\
&
\qquad \qquad
\cdot
\left.\frac{\Phi^{k} (i \sqrt{2(v-w)}) + \Phi^{k} (-i \sqrt{2(v-w)})}{2\sqrt{v-w}}
\right)
\dd w.  \label{eq:F_k_n_minus_k}
\end{align}
\end{thm}

\subsection{Cones, solid angles and intrinsic volumes}\label{sec:geometric_g_n}
The function $g_n$ appeared (in many different parametrisations) in connection with internal and external angles of regular simplices and generalized orthants, but the results are somewhat scattered through the literature and especially the properties of $g_n(r)$ at negative values of $r$ do not seem to be widely known. In the following two sections we shall provide an overview of what is known about $g_n$, state some new results and fix the notation needed for the proof of Theorem~\ref{theo:main1}.

A non-empty subset $C\subset \R^N$ is called a \textit{convex cone} if for every $x,y\in C$ and $\lambda,\mu>0$ we have $\lambda x + \mu y\in C$. In the following we restrict our attention to \textit{polyhedral} cones (or just cones, for short) which are defined as intersections of finitely many closed halfspaces whose boundaries contain the origin. The linear hull of $C$, i.e.\ the smallest linear space containing $C$, is denoted by $L(C)$.
Letting $Z$ be a standard Gaussian random vector on $L(C)$, the \textit{solid angle} of the cone $C$ is defined  as
$$
\alpha(C)=\P[Z\in C].
$$
In fact, the same formula remains true if $Z$ has any rotationally invariant distribution on $L(C)$.
Note that we measure the solid angle w.r.t.\ the linear hull $L(C)$ as the ambient space, so that the solid angle is never $0$, even for cones with empty interior.

Denote the standard scalar product on $\R^N$ by $\langle\cdot,\cdot\rangle$ and let  $e_1,\ldots,e_{N}$ be the standard basis of $\R^{N}$.  Fix any $r\geq -1/n$ and consider $n$ vectors $u_1,\ldots,u_n$ in $\R^N$, $n\leq N$, such that
$$
\langle u_i, u_j\rangle
=
\begin{cases}
1+r, &\text{ if } i=j,\\
r, &\text{ if } i\neq j,
\end{cases}
\quad
i,j\in \{1,\ldots,n\}.
$$
Denote the cone spanned by the vectors $u_1,\ldots,u_n$ by
\begin{equation}\label{eq:C_n_r_def}
C_n(r) := \pos (u_1,\ldots,u_n) := \{\alpha_1 u_1+\ldots+\alpha_n u_n\colon \alpha_1,\ldots,\alpha_n\geq 0\}.
\end{equation}
The specific choice of the vectors $u_1,\ldots,u_n$ as well as the dimension $N$ of the ambient space will be of minor importance because we are interested in the isometry type of the cone only. For $r=0$, the cone $C_n(0)$ is isometric to the positive orthant $\R_+^n$. Vershik and Sporyshev~\cite{vershik_sporyshev} called $C_n(r)$ the contracted (respectively, extended) orthant if $r<0$ (respectively, $r>0$). The extremal cases $r\to \infty$ and $r=-1/n$ correspond to a ray and a half-space, respectively.
\begin{prop}\label{prop:geometric_interpretation_g_n}
For all $r > -1/n$, the solid angle of the cone $C_n(r)$ is given by
$$
\alpha (C_n(r)) = g_n\left(-\frac{r}{1+ n r}\right).
$$
\end{prop}

This fact can be used to relate $g_n(r)$ to the volume of a regular spherical simplex. These volumes have been much studied since~\citet{Schl50}.

\begin{thm}\label{theo:spherical_simplex}
Let $S_n(\ell)$ be a regular spherical simplex, with $n$ vertices and side length $\ell$,  on the unit sphere $\bS^{n-1}$. That is, the geodesic distance between any two vertices of the simplex is $\ell\in (0, \arccos (-\frac 1 {n-1}))$. Then, the spherical volume of $S_n(\ell)$ is given by
$$
\Vol_{n-1}({S_n(\ell)}) = \Vol_{n-1}(\bS^{n-1}) \cdot g_n\left(- \frac{\cos \ell}{1+(n-1)\cos \ell}\right).
$$
More concretely, writing $r_*:= -\frac{\cos \ell}{1+ (n-1)\cos \ell}$, we have
\begin{align}
\Vol_{n-1}({S_n(\ell)}) &= \frac{2}{\Gamma(\frac n2) \sqrt \pi} \int_\infty^{\infty} \left(\int_{-\infty}^{x\sqrt{r_*}}\eee^{-z^2} \dd z\right)^n \eee^{-x^2} \dd x, \quad\text{if } -\frac{1}{n-1} <  \cos \ell \leq 0,\label{eq:vol_spher_simplex1}\\
\Vol_{n-1}({S_n(\ell)}) &=
\frac {2}{\sqrt{\pi}} \int_0^{\infty} \Re \left[\left(\frac 12 + \frac{i}{\sqrt{\pi}} \int_0^{x\sqrt{-r_*} } \eee^{z^2} \dd z\right)^n\right]\eee^{-x^2}\dd x, \quad \text{if } \cos \ell\geq 0. \label{eq:vol_spher_simplex2}
\end{align}
\end{thm}
\begin{proof}
Let $u_1,\ldots, u_n$ be as above. Observe that $u_1/\sqrt{1+r},\ldots,u_n/\sqrt{1+r}$ can be viewed as vertices of $S_n(\arccos \frac{r}{1+r})$. So, choose $r  > -1/n$ such that $\ell = \arccos \frac{r}{1+r}$.  Then, $r=\frac{\cos \ell} {1-\cos \ell}$ and
$$
\frac{\Vol_{n-1}({S_n(\ell)})}{\Vol_{n-1}(\bS^{n-1})} = \alpha(C_n(r)) = g_n\left(-\frac{r}{1+ n r}\right) = g_n\left(- \frac{\cos \ell}{1+(n-1)\cos \ell}\right)
$$
by Proposition~\ref{prop:geometric_interpretation_g_n}.
\end{proof}
Formula~\eqref{eq:vol_spher_simplex1} can be found in the book of \citet[Satz 3 on p.~283]{boehm_hertel_book} or in the works of \citet{ruben, ruben_moments} and~\citet{hadwiger}.  Note that~\cite{boehm_hertel_book} uses a different parametrisation for $S_n(\ell)$; see~\cite[Satz~2 on p.~277]{boehm_hertel_book} for the relation between both parametrisations. Observe also that the Schl\"afli function $\mathbf {f}^{(n)}(\alpha)$ used in~\cite{boehm_hertel_book}  is related to $g_n$ via
$$
\mathbf {f}^{(n)}(\alpha) = \frac {2^n}{n!} g_n\left(-\frac{\cos 2\alpha}{1+\cos 2\alpha}\right)
$$
as one can see by comparing~\cite[Satz~2 on p.~279]{boehm_hertel_book} with Theorem~\ref{theo:spherical_simplex}. The case $\cos \ell >0$ is missing in~\cite{boehm_hertel_book} and in many other references on the subject.
Formula~\eqref{eq:vol_spher_simplex2} was proved by \citet[Corollary~3 on p.~192]{vershik_sporyshev}; see also~\cite{boeroeczky_henk,donoho_tanner1,donoho_tanner2} for asymptotic results.

To proceed, we need to recall some notions related to solid angles. A \textit{polyhedral set} is an intersection of finitely many closed half-spaces (whose boundaries need not pass through the origin).
If a polyhedral set is bounded, it is called a \textit{polytope}. Polyhedral cones are special cases of polyhedral sets.
Denote by $\cF_k(P)$ the set of $k$-dimensional faces of a polyhedral set $P$. The \textit{tangent cone} at a face $F\in \cF_k(P)$ is defined by
\begin{equation}\label{eq:def_tangent_cone}
T_F = T_F(P) = \{v\in\R^n\colon f_0 +\eps v \in P \text{ for some } \eps>0\},
\end{equation}
where $f_0$ is any point in the relative interior of $F$, i.e.\ the interior of $F$ taken w.r.t.\ its affine hull. The \textit{normal cone} at the face $F\in \cF_k(P)$ is defined as the polar of the tangent one, that is
\begin{equation}\label{eq:def_normal_cone}
N_F = N_F(P) = T_F^\circ (P) = \{w\in\R^n \colon \langle w, u\rangle\leq 0 \text{ for all } u\in T_F(P)\}.
\end{equation}
For certain special values of $r$ one can interpret $g_n(r)$ as inner or normal solid angles at the faces of the regular simplex. The \textit{inner} and \textit{normal} (or \textit{external}) solid angles of $P$ at $F$ are defined as the solid angles of the cones $T_F(P)$ and $N_F(P)$, respectively.
\begin{prop}\label{prop:angles_simplex}
Let $\Delta_n:= \conv(e_1,\ldots,e_n)$ be the $(n-1)$-dimensional regular simplex in $\R^n$.
\begin{itemize}
\item[(a)] The normal solid angle at any $(k-1)$-dimensional face of $\Delta_n$ equals
$$
g_{n-k}\left(\frac 1k\right) = \frac {1}{\sqrt{\pi}} \int_{-\infty}^{\infty} \left(\frac{1}{\sqrt{\pi}} \int_{-\infty}^{x/\sqrt{k}} \eee^{-z^2} \dd z\right)^{n-k}\eee^{-x^2}\dd x.
$$
\item[(b)] The inner solid angle at any $(k-1)$-dimensional face of $\Delta_n$ equals
$$
g_{n-k}\left(-\frac 1n\right) =
\frac {2}{\sqrt{\pi}} \int_0^{\infty} \Re \left[\left(\frac 12 + \frac{i}{\sqrt{\pi}} \int_0^{x/\sqrt{n}} \eee^{z^2} \dd z\right)^{n-k}\right]\eee^{-x^2}\dd x.
$$
\end{itemize}
\end{prop}
Both parts were known; see~\cite{hadwiger} and~\cite{ruben} for part~(a) and~\cite[Section~4]{rogers} (where the method used was attributed to H.\ E.\ Daniels) as well as~\cite[Lemma~4]{vershik_sporyshev} for part~(b). A formula for the normal solid angles of crosspolytopes (which is similar to part~(a)) was derived in~\cite{betke_henk}.

The next proposition provides a geometric interpretation of $b_{n,k}(\sigma^2)$. The $k$-th \textit{conic intrinsic volume} of a cone $C$ is given by
$$
\upsilon_k(C) = \sum_{F\in \cF_k(C)} \alpha(F) \alpha(N_F(C)),
$$
where we recall that $\alpha(F)$ is the solid angle of the cone $F$ measured w.r.t.\ the linear hull of $F$. See~\cite{amelunxen_lotz} for equivalent definitions and properties.
\begin{prop}\label{prop:upsilon_n_k}
For every $r >  -1/n$ and $k\in \{0,\ldots,n\}$, the $k$-th conic intrinsic volume of the cone $C_n(r)$ is given by
$$
\upsilon_k(C_n(r)) = b_{n,k}(r) = \binom nk g_k\left(-\frac{r}{1+kr}\right) g_{n-k} \left(\frac{r}{1+kr}\right).
$$
\end{prop}
\begin{rem}\label{rem:conic_gauss_bonnet}
As a consequence of the Gauss--Bonnet formula for conic intrinsic volumes, see~\cite[Theorem~6.5.5]{SW08} or~\cite[Corollary~4.4]{amelunxen_lotz}, we obtain the identities
$$
\sum_{k=0}^{\lfloor n/2 \rfloor} b_{n,2k}(r) = \sum_{k=0}^{\lfloor (n-1)/2\rfloor} b_{n,2k+1}(r) = \frac 12.
$$
In particular, the numbers $b_{n,0}(r), \ldots, b_{n,n}(r)$ define a probability distribution on $\{0,\ldots,n\}$, a fact which is not evident in view of the expression for $g_n$ given in~\eqref{eq:g_n_first_formula1} and~\eqref{eq:g_n_first_formula2}. For $r=0$ (in which case $C_n(0)$ is the positive orthant $\R^n_+$) this distribution reduces to the binomial one with parameters $(n, 1/2)$ because $g_n(0)= 2^{-n}$ by Proposition~\ref{prop:g_n_properties}~(d), below.
\end{rem}

\subsection{Properties of \texorpdfstring{$g_n$}{gn}}\label{sec:explicit_g_n}
Next we give a formula for $g_n(r)$ which may be more convenient than its definition~\eqref{eq:def_g}.
Recall that $\Phi$ denotes the standard normal distribution function, see~\eqref{eq:Phi_def}, viewed as an analytic function on the entire complex plane.

\begin{prop}\label{prop:g_n_properties}
The function $g_n: [-\frac 1n,\infty) \to [0,1]$ defined in~\eqref{eq:def_g} has the following properties.
\begin{itemize}
\item[(a)]
For all $n\in\boldN$ and $r\geq -\frac 1n$,
\begin{align}
g_n(r)
&=
\frac 1 {\sqrt {2\pi}} \int_{-\infty}^{\infty} \Phi^n (\sqrt r x) \eee^{-x^2/2} \dd x\notag \\
&=
\frac 1 {\sqrt {2\pi}} \int_{0}^{\infty} \left(\Phi^n (\sqrt r x) + \Phi^n (-\sqrt r x)\right) \eee^{-x^2/2} \dd x,\label{eq:g_n_r_part_a}
\end{align}
where, in the case of negative $r$, we use the convention $\sqrt{r} = i \sqrt{-r}$. In fact, the right-most expression in~\eqref{eq:g_n_r_part_a} defines $g_n$ as an analytic function on the half-plane $\Re r>-\frac 1n$.
\item[(b)]
For all $n\in \{2,3,\ldots\}$ and $r> -\frac 1n$ we have
$$
g_n'(r) = \frac{n(n-1)}{4\pi (r+1)\sqrt{2r+1}} g_{n-2} \left(\frac{r}{2r+1}\right).
$$
\item[(c)]
$g_0(r) = 1$  (by definition) and $g_1(r) = \frac 12$. 
\item[(d)] For every $n\in\boldN$, we have
$$
g_n(0) = 2^{-n},
\;\;\;
g_n(1) = \frac 1 {n+1},
\;\;\;
\lim_{r\to+\infty}g_n(r)= \frac 12.
$$
\item[(e)]
For $n\in \{2,3,\ldots\}$ we have $g_n(-\frac 1n) = 0$.
\item[(f)]
$g_2(r) = \frac 14 + \frac 1{2\pi} \arcsin \frac{r}{1+r}$ and
$g_3(r) = \frac 18 + \frac 3{4\pi} \arcsin \frac{r}{1+r}$.

\item[(g)] For every fixed $n\in\boldN$ we have
$$
g_n\left( -\frac 1n + \eps\right) \sim  \frac{n^n \Gamma(n/2)}{2 \,  \pi^{n/2} \Gamma(n) \sqrt n} \cdot \eps^{(n-1)/2}, \quad \eps \downarrow 0.
$$
\item[(h)]
For all $n\in\boldN$, the functions $g_{2n}$ and $g_{2n+1}$ admit extensions to analytic functions on some unramified cover of $\C\bsl \{-1/k\colon  k=1,\ldots, 2n\}$.
\end{itemize}
\end{prop}
\begin{rem}
Special values of $g_n$ listed in Parts (d) and (e) were known to Schl\"afli~\cite[p.~267]{Schl50}; see also~\cite[pp.~285--286]{boehm_hertel_book}.
Part (b) is a consequence of the Schl\"afli differential relation; see~\citet[Satz~2, p.~279]{boehm_hertel_book}. For completeness, we shall provide a self-contained proof of Proposition~\ref{prop:g_n_properties} in Section~\ref{subsec:g_n_properties}.
\end{rem}

\begin{rem}
Using the fact that $\overline{\Phi(iz)} = \Phi(-iz)$ for $z\in\R$, one can state~\eqref{eq:g_n_r_part_a} in the case of real $r \in [-\frac 1n, 0]$ as follows:
\begin{align*}
g_n(r)
&=
\frac 2 {\sqrt {2\pi}}  \int_{0}^{\infty} \Re(\Phi^n ( i \sqrt {-r} x)) \eee^{-x^2/2} \dd x\\
&=
\frac {2}{\sqrt{\pi}} \int_0^{\infty} \Re \left[\left(\frac 12 + \frac{i}{\sqrt{\pi}} \int_0^{x\sqrt{-r}} \eee^{z^2} \dd z\right)^n\right]\eee^{-x^2}\dd x,  
\end{align*}
a formula obtained by \citet[Corollary 3 on p. 192]{vershik_sporyshev}.
\end{rem}
\begin{rem}
Taking $r=-1/n$ in~\eqref{eq:g_n_r_part_a}, making the change of variables $y=x/\sqrt n$ and using that $g_n(-1/n) = 0$ for $n\geq 2$, we obtain the curious identity
\begin{equation}\label{eq:identity_integral_phi}
\int_{-\infty}^{+\infty} \left(\Phi(iy)\eee^{-y^2/2}\right)^n \dd y = 0,
\quad
n = 2,3,\ldots.
\end{equation}
Using induction and partial integration we shall extend this as follows.
\end{rem}
\begin{prop}\label{prop:integral_phi}
For all $m\in\boldN_0$ and all $n = m+2,m+3,\ldots$ we have
\begin{equation}\label{eq:int_erfi}
\int_{-\infty}^{+\infty} y^{m} \left(\Phi(iy)\eee^{-y^2/2}\right)^n \dd y = 0.
\end{equation}
Also, for all $m\in\boldN_0$  we have
\begin{equation}\label{eq:int_erfi_cauchy}
\int_{0}^{+\infty} y^{m} \left(\left(\Phi(iy)\eee^{-y^2/2}\right)^{m+1} + (-1)^m \left(\Phi(-iy)\eee^{-y^2/2}\right)^{m+1}\right) \dd y
=
\sqrt{\frac \pi 2} \left(\frac i {\sqrt{2\pi}}\right)^m.
\end{equation}
\end{prop}
For $n\leq m+1$ the integral in~\eqref{eq:int_erfi} diverges  since $\Phi(iy) \sim \frac1 {\sqrt{2\pi} iy}\eee^{y^2/2}$ as $y\to\infty$, $y\in\R$; see~\cite[Eq. 7.1.23 on p. 298]{abramowitz_stegun}.  Equation~\eqref{eq:int_erfi_cauchy} states a formula for the Cauchy principal value which is well defined  for $n=m+1$.

\subsection{Expected number of \texorpdfstring{$k$}{k}-faces}
Let $f_k(\cP_{n,d})$ be the number of $k$-dimensional faces of the Gaussian polytope $\cP_{n,d}$. Recall the notation $p_{n,d} (\sigma^2) = \P[\sigma X \notin \cP_{n,d}]$, where $X$ is a standard normal vector in $\R^d$ independent of $\cP_{n,d}$.
With the aid of the Blaschke--Petkantschin formula, \citet[Theorem~3.2]{HMR04} showed  that
\begin{equation}\label{eq:E_f_k_hug_etal}
\boldE f_k ( \cP_{n,d}) = \binom{n}{k+1} p_{n-k-1, d-k}\left(\frac 1 {k+1}\right).
\end{equation}
Using this formula, they proved an asymptotic result of the form
\begin{equation}\label{eq:asympt_E_f_k_n}
\boldE f_k ( \cP_{n,d}) \simn  \bar c_{(k,d)} (\log n)^{(d-1)/2}.
\end{equation}
where $\bar c_{(k,d)}$ is an explicit constant only depending on $d$ and $k$. With the aid of Theorem~\ref{theo:main1} one can derive the following explicit formula.
\begin{thm}\label{theo:expected_faces}
For every $k=0,\ldots, n-1$ we have
\begin{equation}\label{eq:E_f_k_P_n_d}
\boldE f_k ( \cP_{n,d}) = \frac{2\cdot n!}{(k+1)!} \sum_{\substack{j= d - 2i > k\\i\in \boldN_0}} \frac{g_{j-1-k}\left(-\frac 1j\right)g_{n-j}\left(\frac 1j\right)}{(j-1-k)!(n-j)!}.
\end{equation}
\end{thm}
\begin{proof}
Combine Theorem~\ref{theo:main1} with~\eqref{eq:E_f_k_hug_etal}.
\end{proof}
\begin{rem}
The quantities $g_{j-1-k}(-\frac 1j)$ and $g_{n-j}(\frac 1j)$ appearing on the right-hand side of~\eqref{eq:E_f_k_P_n_d} are certain inner/normal solid angles of regular simplices; see Proposition~\ref{prop:angles_simplex}. With this interpretation, formula~\eqref{eq:E_f_k_P_n_d} is due to \citet{AS92}.
\end{rem}
\begin{rem}
More generally, \citet{HMR04} considered also the functional
$$
T^{d,k}_{0,b} (\cP_{n,d}) = \sum_{F\in \cF_k(\cP_{n,d})} (\Vol_k(F))^b,  \quad b\geq 0,
$$
which reduces to $f_{k} (\cP_{n,d})$ for $b=0$. For $k=d-1$ and $b=1$ one gets the surface area.   \citet{HMR04} showed that
$$
\boldE T^{d,k}_{0,b} (\cP_{n,d}) = \boldE f_k ( \cP_{n,d})\cdot  \left(\frac{\sqrt {k+1}}{k!}\right)^b \prod_{j=1}^k \frac{\Gamma\left(\frac{d+b+1-j}{2}\right)}{\Gamma\left(\frac{d+1-j}{2}\right)}
.
$$
Thus, an explicit formula for $\boldE T^{d,k}_{0,b} (\cP_{n,d})$ can be obtained from Theorem~\ref{theo:expected_faces}.
\end{rem}


The following  fixed $r$ asymptotics for $g_n(r)=\frac 1 {\sqrt {2\pi}} \int_{-\infty}^{\infty} \Phi^n (\sqrt r x) \eee^{-x^2/2} \dd x$ was derived in~\cite[pp.~44-45]{raynaud} and~\cite[Lemma~5]{vershik_sporyshev_asymptotic}.
\begin{prop}\label{prop:asymptotics_g}
For any fixed $r>0$ we have
$$
g_n(r) \simn \Gamma(1/r) r^{-1/2} n^{-1/r} (4\pi\log n)^{\frac 1{2r} - \frac 12}.
$$
\end{prop}
This can be used to compute the large $n$ asymptotics of the probability that $\sigma X \notin \cP_{n,d}$.
\begin{cor}\label{cor:asympt_p_n_d}
Fix any $\sigma^2>0$ and write $r=\frac{\sigma^2}{1+(d-1)\sigma^2}$.  Then,
$$
p_{n,d}(\sigma^2) = \P[\sigma X \notin \cP_{n,d}] \simn \frac {n^{-1/\sigma^2}} {(d-1)!}  g_{d-1}(-r)
\Gamma\left(\frac 1 r\right) r^{-1/2}
(4\pi\log n)^{\frac{1}{2r} - \frac 12}.
$$
\end{cor}
\begin{proof}
It follows from Proposition~\ref{prop:asymptotics_g} and~\eqref{eq:def_b_k} that  for every fixed $k\in\boldN_0$ and $\sigma^2>0$,
$$
b_{n,k} (\sigma^2)
\simn   \Gamma\left(k + \sigma^{-2}\right) \frac {\sqrt{k + \sigma^{-2}}} {k!} g_k\left(-\frac{\sigma^2}{1+k\sigma^2}\right) n^{-\frac 1{\sigma^2}}
(4\pi\log n)^{\frac{1}{2\sigma^2} + \frac {k-1}{2}}.
$$
In particular, $b_{n,d-m} (\sigma^2) = o(b_{n, d-1}(\sigma^2))$ for all $m\geq 2$ and as $n\to\infty$. Theorem~\ref{theo:main1} yields
$p_{n,d}(\sigma^2) \sim 2 b_{n,d-1} (\sigma^2)$, from which the required asymptotics follows.
\end{proof}
\begin{rem}
Applying Corollary~\ref{cor:asympt_p_n_d} to the right-hand side of~\eqref{eq:E_f_k_hug_etal} we deduce the asymptotic formula
\begin{align*}
\boldE f_k ( \cP_{n,d})
&=
\binom{n}{k+1} p_{n-k-1, d-k}\left(\frac 1 {k+1}\right)\\
&\simn \frac{2}{\sqrt d} \binom d{k+1} g_{d-1-k} \left(-\frac 1d\right) (4\pi \log n)^{(d-1)/2},
\end{align*}
which recovers a result of \citet{affentranger} (see also~\cite{AS92,baryshnikov_vitale} and~\cite{HMR04}) stated in~\eqref{cor:asympt_p_n_d}.
\end{rem}

\subsection{Expected volume}
Let us derive from our results the following  formula for the expected volume of the Gaussian polytope due to~\citet{efron}:
\begin{equation}\label{eq:efron}
\boldE \Vol_d (\cP_{n,d}) =  \frac {\pi^{\frac d2}}{\Gamma(\frac d2 +1)}  \cdot \frac{n!}{d! (n-d-1)!} \int_{-\infty}^{\infty} \Phi^{n-d-1}(t) \varphi^{d+1} (t) \dd t.
\end{equation}
In fact, \citet{efron} proved the formula for $d=2$ and stated it for general $d$; another proof (valid for arbitrary $d$) can be found in~\cite{kabluchko_zaporozhets_intrinsic}.

Since the surface measure of the unit sphere in $\R^d$ equals $\omega_d = 2\pi^{d/2}/\Gamma(d/2)$, we can express the expected volume of the Gaussian polytope as follows:
\begin{align*}
\boldE \Vol_d (\cP_{n,d})
&=
\int_0^\infty (1-f_{n,d} (r))  \frac{2\pi^{d/2}}{\Gamma(d/2)} r^{d-1} \dd r\\
&=
\frac {(2\pi)^{d/2}}{\Gamma(d/2)}\int_0^\infty (1-f_{n,d} (\sqrt{2u}))  u^{(d/2)-1}\dd u,
\end{align*}
where we have made the change of variables $r= \sqrt{2u}$. On the other hand, by Corollary~\ref{cor:intensity} together with the identity $\sum_{k=0}^n b_{n,k}(1/\lambda) =1$ (see Remark~\ref{rem:conic_gauss_bonnet}), we can write
$$
\int_0^\infty (1 - f_{n,d}(\sqrt{2u})) u^{(d/2) - 1} \eee^{-\lambda u} \dd u  =  2 \Gamma (d/2) \lambda^{-d/2}  (b_{n,d+1} (1/\lambda) + b_{n,d+3} (1/\lambda) + \ldots)
$$
for all $\lambda>0$. Hence, by the monotone convergence theorem,
\begin{equation}\label{eq:E_Vol}
\boldE \Vol_d (\cP_{n,d}) = 2^{(d/2) + 1} \pi^{d/2} \lim_{\lambda \downarrow 0} \lambda^{-d/2}  (b_{n,d+1} (1/\lambda) + b_{n,d+3} (1/\lambda) + \ldots).
\end{equation}
Recall from~\eqref{eq:def_b_k} that for every fixed $k\in \{0,\ldots,n\}$,
$$
b_{n,k}(1/\lambda) = \binom nk g_k\left(-\frac{1}{\lambda+k}\right) g_{n-k} \left(\frac{1}{\lambda+k}\right).
$$
Clearly,
$$
\lim_{\lambda \downarrow 0} g_{n-k} \left(\frac{1}{\lambda+k}\right) = g_{n-k} \left(\frac 1k \right),
$$
whereas Proposition~\ref{prop:g_n_properties} (g) yields
$$
g_k\left(-\frac{1}{\lambda+k}\right)  = g_k\left(-\frac{1}{k} + \frac{\lambda(1 + o(1))}{k^2}\right) \sim \frac {\Gamma(k/2) \sqrt k}{2 \pi^{k/2} \Gamma(k)} \lambda^{(k-1)/2} \text{ as } \lambda \downarrow 0.
$$
Taking everything together, we obtain
$$
b_{n,k}(1/\lambda) \sim  \binom nk   \frac {\Gamma(k/2) \sqrt k}{2 \pi^{k/2} \Gamma(k)} g_{n-k} \left(\frac 1k \right) \lambda^{(k-1)/2} \text{ as } \lambda\downarrow 0.
$$
So, in the sum on the right-hand side of~\eqref{eq:E_Vol} the term $b_{n,d+1}(1/\lambda)\sim \text{const} \cdot \lambda^{d/2}$ dominates and we obtain
\begin{align*}
\boldE \Vol_d (\cP_{n,d})
&=
(2\pi)^{d/2} \binom{n}{d+1} \frac{\Gamma(\frac{d+1}{2}) \sqrt{d+1}}{ \pi^{\frac {d+1}{2}} \Gamma(d+1)} g_{n-d-1}\left(\frac 1 {d+1}\right)\\
&=
 \binom{n}{d+1} \frac{\sqrt{d+1}}{ 2^{\frac d2} \Gamma(\frac d2 + 1)} g_{n-d-1}\left(\frac 1 {d+1}\right),
\end{align*}
where we used the Legendre duplication formula.  Recalling formula~\eqref{eq:g_n_r_part_a} for $g_n(r)$, and performing some simple transformations, we arrive at~\eqref{eq:efron}. Proposition~\ref{prop:asymptotics_g} yields the following asymptotic formula due to~\citet[Theorem 4]{affentranger}:
$$
\boldE \Vol_d (\cP_{n,d}) \sim \frac{\pi^{d/2}}{\Gamma(\frac d2 +1)} (2\log n)^{d/2}
$$
as $n\to\infty$, while $d$ stays fixed. A more refined asymptotics was derived in~\cite{kabluchko_zaporozhets_intrinsic}.




\subsection{Low-dimensional examples} Let $d=2$. Theorem~\ref{theo:main1} simplifies to
\begin{align*}
\P[\sigma X \notin \conv (X_1,\ldots,X_n)]
&=
2b_{n,1}(\sigma^2) = 2 n g_1\left(-\frac{\sigma^2}{1+\sigma^2}\right) g_{n-1}\left(\frac{\sigma^2}{1+\sigma^2}\right)\\
&=
ng_{n-1}\left(\frac{\sigma^2}{1+\sigma^2}\right).
\end{align*}
We obtain the formula
$$
\P[\sigma X \notin \conv (X_1,\ldots,X_n)]
=
\frac{n}{\sqrt{2\pi}}\int_{-\infty}^{+\infty} \Phi^{n-1} \left(\frac{\sigma x}{\sqrt{1+\sigma^2}}\right) \eee^{-x^2/2} \dd x.
$$
The next theorem gives an explicit formula for the non-absorption probability in the case $d=2$.
\begin{thm}\label{theo:intensity_d_2}
Let $\xi,\xi_1,\ldots,\xi_n$ be standard normal random variables and let $W$ be a random variable with the arcsine density $\frac{1}{\pi} (1-y^2)^{-1/2} \dd y$ on $[-1,1]$, all variables being independent. Define $M_{k}= \max\{\xi_1,\ldots,\xi_{k}\}$.  Then, for all $u\geq 0$,
\begin{align*}
f_{n,2}(\sqrt{2u})
&= \P[ M_{n}^2 + \xi^2 \leq 2u] + \frac {\dd}{\dd u}\P[ M_{n}^2 + \xi^2 \leq 2u]\\
&=  \P[ M_{n}^2 + \xi^2 \leq 2u] + n\eee^{-u} \P[M_{n-1} \leq \sqrt{2u} W].
\end{align*}
That is, $f_{n,2}(\sqrt{2u})$ is the sum of the distribution function and the density of the random variable $\frac 12 (M_n^2+\xi^2)$ at $u$.
\end{thm}
Two-dimensional absorption probabilities were studied in a very general setting by Jewell and Romano~\cite{jewell_romano,jewell_romano_inclusion}, but their method does not seem to yield an explicit formula like that in Theorem~\ref{theo:intensity_d_2}.

Consider now the case $d=3$. Using first Theorem~\ref{theo:main1}, \eqref{eq:def_b_k}, and then Proposition~\ref{prop:g_n_properties}, we arrive at
\begin{align*}
\lefteqn{\P[\sigma X \notin \conv (X_1,\ldots,X_n)]}\\
&= 2b_{n,2}(\sigma^2) + 2b_{n,0}(\sigma^2)\\
&=
n(n-1)  g_2\left(-\frac{\sigma^2}{1+2\sigma^2}\right) g_{n-2}\left(\frac{\sigma^2}{1+2\sigma^2}\right)
+ 2 g_n(\sigma^2)\\
&=
\frac{n(n-1)}{\sqrt{2\pi}}
\left(\frac 14 - \frac 1{2\pi} \arcsin \frac{\sigma^2}{1+\sigma^2}\right)
\int_{-\infty}^{+\infty} \Phi^{n-2} \left(\frac{\sigma x}{\sqrt{1+2\sigma^2}}\right) \eee^{-x^2/2} \dd x
\\& + \frac {2}{\sqrt{2\pi}} \int_{-\infty}^{+\infty} \Phi^n (\sigma x) \eee^{-x^2/2}\dd x,
\end{align*}
but we were unable to invert the Laplace transform to obtain a formula for $f_{n,3}$ similar to that of Theorem~\ref{theo:intensity_d_2}. Similarly, for $d=4$ we obtain
\begin{align*}
\lefteqn{\P[\sigma X \notin \conv (X_1,\ldots,X_n)]}\\
&= 2b_{n,3}(\sigma^2) + 2b_{n,1}(\sigma^2)\\
&=
\frac{n(n-1)(n-2)}{3}  g_3\left(-\frac{\sigma^2}{1+3\sigma^2}\right) g_{n-3}\left(\frac{\sigma^2}{1+3\sigma^2}\right) + ng_{n-1}\left(\frac{\sigma^2}{1+\sigma^2}\right)\\
&=
\frac{n(n-1)(n-2)}{3\sqrt{2\pi}}
\left(\frac 18 - \frac 3{4\pi} \arcsin \frac{\sigma^2}{1+2\sigma^2}\right)
\int_{-\infty}^{+\infty} \Phi^{n-3} \left(\frac{\sigma x}{\sqrt{1+3\sigma^2}}\right) \eee^{-x^2/2} \dd x
\\&\phantom{=}+ \frac{n}{\sqrt{2\pi}}\int_{-\infty}^{+\infty} \Phi^{n-1} \left(\frac{\sigma x}{\sqrt{1+\sigma^2}}\right) \eee^{-x^2/2} \dd x,
\end{align*}


Taking $\sigma=1$ in Theorem~\ref{theo:main1} yields the \textit{probability content} of the Gaussian polytope which is defined as
$$
C_{n,d} := \P[X \in \conv (X_1,\ldots,X_n)] = 1 - \P[X\notin \conv (X_1,\ldots,X_n)].
$$
For $d=2,3,4$ we obtain the formulas
\begin{align*}
C_{n,2} &= 1- \frac{n}{\sqrt{2\pi}}\int_{-\infty}^{+\infty} \Phi^{n-1} \left(\frac{x}{\sqrt{2}}\right) \eee^{-x^2/2} \dd x,\\
C_{n,3} &= 1-
\frac{n(n-1)}{6\sqrt{2\pi}}
\int_{-\infty}^{+\infty} \Phi^{n-2} \left(\frac{x}{\sqrt 2}\right) \eee^{-x^2/2} \dd x - \frac {2}{n+1},\\
C_{n,4} &= C_{n,2} -
\frac{n(n-1)(n-2)}{3\sqrt{2\pi}} \left( \frac 18 - \frac 3{4\pi} \arcsin \frac 13\right)
\int_{-\infty}^{+\infty} \Phi^{n-3} \left(\frac{x}{2}\right) \eee^{-x^2/2} \dd x.
\end{align*}
The formulas for $d=2,3$ were obtained by~\citet{efron}, Equations (7.5) and (7.6) on p.~341.

\subsection{Absorption probability in the Goodman--Pollack model}
Let $v_1,\ldots,v_{n+1}$ be the vertices of an $n$-dimensional regular simplex inscribed into the unit sphere $\mathbb S^{n-1} \subset \R^{n}$. That is,  $|v_i|=1$ for all $1\leq i\leq n+1$ and  $\rho := \langle v_i,v_j\rangle = -1/n$ for all $1\leq i < j \leq n+1$.  Let $O$ be a random orthogonal matrix sampled according to the Haar probability measure on the orthogonal group $O(n)$. Consider the randomly rotated regular simplex with vertices $Ov_1,\ldots,Ov_{n+1}$ and project it onto some fixed $d$-dimensional linear subspace $V_d\subset \R^{n}$. The choice of $V_d$ is irrelevant, so that we shall assume that $V_d\equiv \R^d$ is spanned by the first $d$ vectors of the standard orthonormal basis $e_1,\ldots,e_n$ of $\R^{n}$.  Denote orthogonal projection onto $V_d$ by $\Pi$. The resulting random polytope
$$
\cQ_{n+1,d} := \conv (\Pi O v_1,\ldots, \Pi O v_{n+1})\subset \R^d
$$
is said to be distributed according to the \textit{Goodman--Pollack model}. \citet{AS92} and \citet{baryshnikov_vitale} observed  that the Gaussian polytope $\cP_{n+1,d}$  is closely related to the Goodman--Pollack polytope $\cQ_{n+1,d}$.  In particular, \citet{baryshnikov_vitale} showed that all functionals which remain invariant under affine transformations of the polytope  (like the number of $k$-faces) have the same distribution in both models. We are interested in the non-absorption probability in the Goodman--Pollack model, that is
$$
f_{n,d}^* (|x|) := \P[x\notin \cQ_{n,d}], \quad x\in\R^d.
$$
Clearly, this functional is not invariant w.r.t.\ the affine transformations of $\cQ_{n,d}$. We cannot compute the non-absorption probability explicitly, but it is possible to evaluate certain integral transform of $f_{n,d}^*$.
\begin{thm}\label{theo:goodman_pollack}
For every $\sigma^2\geq 0$ we have
$$
\frac 1 {B\left(\frac d2, \frac{n-d}{2}\right)} \int_{0}^{\infty} \frac{u^{d-1}}{(1+u^2)^{n/2}}  f_{n,d}^* \left(\sqrt{\frac{n\sigma^2 + 1}{n-1}} \cdot u\right) \dd u = b_{n,d-1}(\sigma^2) + b_{n,d-3}(\sigma^2) +\ldots.
$$
\end{thm}
Here,  $B(\cdot, \cdot)$ denotes the Euler  Beta function. The proof of Theorem~\ref{theo:goodman_pollack} will be given in Section~\ref{sec:proof_goodman_pollack}.
\begin{rem}\label{rem:goodman_pollack_standard_simplex}
It is also possible to consider random projections $\tilde \cQ_{n,d}$ of the random orthogonal transformation of the regular simplex $\conv (e_1,\ldots,e_n)$ inscribed into the unit sphere $\mathbb {S}^{n-1} \subset \R^n$. For the non-absorption probability $\tilde f_{n,d}^*(|x|):= \P[x\notin \tilde \cQ_{n,d}]$, $x\in\R^d$, one can obtain
$$
\frac 1 {B\left(\frac d2, \frac{n-d+1}{2}\right)} \int_{0}^{\infty} \frac{u^{d-1} \tilde f_{n,d}^* (\sigma u)}{(1+u^2)^{(n+1)/2}}   \dd u = b_{n,d-1}(\sigma^2) + b_{n,d-3}(\sigma^2) +\ldots
$$
by a slight simplification of the proof of Theorem~\ref{theo:goodman_pollack}; see Remark~\ref{rem:goodman_pollack_proof}.
\end{rem}

\section{Proof of Theorem~\ref{theo:main1}} \label{sec:proof_main}


\subsection{Reduction to intrinsic volumes}
We can replace $X$ by $-X$ because by the symmetry of the Gaussian distribution
$$
\P[\sigma X \notin \conv (X_1,\ldots, X_n)] = \P[-\sigma X \notin \conv (X_1,\ldots, X_n)].
$$
Clearly, $-\sigma X \notin \conv (X_1,\ldots,X_n)$ if and only if $0\notin \conv (X_1+\sigma X, \ldots, X_n+\sigma X)$. This, in turn, is equivalent to the following condition:
$$
\alpha_1 X_1 + \ldots + \alpha_n X_n + (\alpha_1+\ldots+\alpha_n) \sigma X = 0,
\;\;
\alpha_1,\ldots,\alpha_n\geq 0
\;\;
\Longrightarrow
\;\;
\alpha_1=\ldots=\alpha_n = 0.
$$
To interpret this geometrically, we consider the following convex cone in the space $\R^{n+1}$:
\begin{equation}\label{eq:C_def}
C:= \{(\alpha_1,\ldots,\alpha_n, (\alpha_1+\ldots+\alpha_n) \sigma ) \colon \alpha_1,\ldots,\alpha_n \geq 0\}.
\end{equation}
This cone is spanned by $e_1+\sigma e_{n+1},\ldots,e_n+\sigma e_{n+1}$ (where $e_1,\ldots,e_{n+1}$ is the standard basis of $\R^{n+1}$) and is therefore isometric to the cone $C_n(\sigma^2)$ introduced in~\eqref{eq:C_n_r_def}.
Let also $U$ be a random linear subspace of $\R^{n+1}$ given by
$$
U=\{(y_1,\ldots,y_{n+1})\in \R^{n+1}\colon y_1X_1+\ldots+y_n X_n + y_{n+1} X = 0\}.
$$
Observe that since $X_1,\ldots,X_n, X$ are i.i.d.\ standard Gaussian random vectors on $\R^d$, where $d < n$,  the linear space $U$ has a.s.\ codimension $d$ and is uniformly distributed on the corresponding linear Grassmannian.
The above discussion shows that
\begin{equation}\label{eq:absorpt_as_intersect}
\P[\sigma X \notin \conv (X_1,\ldots,X_n)]
= \P[U\cap C = \{0\}].
\end{equation}
The next result, known as the conic Crofton formula~\cite[pp.~261--262]{SW08} or~\cite[Corollary~5.2]{amelunxen_lotz}, is of major importance for us.
\begin{thm}\label{theo:crofton_conic}
Let $C\subset \R^N$ be a convex polyhedral cone which is not a linear subspace. If $U\subset \R^N$ is a uniformly distributed linear subspace of codimension $d$, then
$$
\P[C\cap U = \{0\}] = 2 (\upsilon_{d-1}(C) + \upsilon_{d-3}(C) + \ldots),
$$
where $\upsilon_0(C), \ldots, \upsilon_N(C)$ are the  conic intrinsic volumes of $C$ given by
\begin{equation}\label{eq:upsilon}
\upsilon_k(C) = \sum_{F\in \cF_k(C)} \alpha(F) \alpha(N_F(C)).
\end{equation}
\end{thm}
Combining~\eqref{eq:absorpt_as_intersect} with the conic Crofton formula, we obtain
$$
\P[\sigma X \notin \conv (X_1,\ldots,X_n)]
= \P[U\cap C = \{0\}]
= 2 (\upsilon_{d-1}(C) + \upsilon_{d-3}(C) + \ldots).
$$

In the following, we shall show that the number of $k$-faces of $C=C_n(\sigma^2)$ is $\binom nk$, and for every $k$-face $F\in \cF_k(C)$ we have
\begin{equation}\label{eq:angles_claim}
\alpha(F) = g_k\left(-\frac{\sigma^2}{1+k\sigma^2}\right),
\quad
\alpha(N_F(C)) = g_{n-k}\left(\frac{\sigma^2}{1+k\sigma^2}\right),
\end{equation}
where $g_k(r)$ is as in Section~\ref{subsec:main_results}. This would prove Theorem~\ref{theo:main1}.

\subsection{The polar cone}
The \textit{polar cone} of  a convex cone $D\subset \R^N$ is defined by
$$
D^\circ = \{x\in\R^N \colon \langle x,y \rangle\leq 0 \text{ for all } y\in D\}.
$$
\begin{prop}\label{prop:dual_cone}
Let $r >-1/n$. The polar cone of $C_n(r)$ taken with respect to the ambient space $L(C_n(r))$ is isometric to $C_{n}(-\frac{r}{1+n r})$. That is to say, there is an orthogonal transformation $O:\R^N\to\R^N$ such that
$$
O(C_n^\circ(r) \cap L(C_n(r))) = C_{n}\left(-\frac{r}{1+n r}\right).
$$
\end{prop}
\begin{proof}
Since $D^{\circ \circ} = D$ and since the transformation $r\mapsto -\frac {r}{1+nr}$ is an involution, it suffices to prove the proposition for $r=\sigma^2 \geq 0$.   Since we work in the linear hull of the cone, there is no restriction of generality in assuming that it has the form $C=C_n(\sigma^2)$ given in~\eqref{eq:C_def}. Thus, $C$ is spanned by the vectors
 $u_1,\ldots,u_n$ given by $u_i = e_i + \sigma e_{n+1}\in \R^{n+1}$, $1\leq i \leq n$. The linear space spanned by $u_1,\ldots,u_n$ is
$$
L(C) =
\{(\alpha_1,\ldots,\alpha_n,\alpha_{n+1}) \in \R^{n+1} \colon \alpha_{n+1} = \sigma(\alpha_1+\ldots + \alpha_{n})\}.
$$
The polar cone of $C=C_n(\sigma^2)$ taken with respect to $L(C)$ as the ambient space is
$$
C^{\circ} \cap L(C) = \{(\alpha_1,\ldots,\alpha_{n+1})\in L(C)\colon \alpha_1+ \sigma \alpha_{n+1} \leq 0,\ldots,\alpha_n + \sigma \alpha_{n+1} \leq 0 \}.
$$
The lineality space of a cone $D$ is defined as $D\cap (-D)$. The lineality space of $C^{\circ} \cap L(C)$ is trivial, namely
$$
\{(\alpha_1,\ldots,\alpha_{n+1})\in L(C)\colon \alpha_1+ \sigma \alpha_{n+1}=0,\ldots,\alpha_n + \sigma \alpha_{n+1} =0 \} = \{0\}.
$$
It follows that the cone $C^{\circ} \cap L(C)$  is spanned by its one-dimensional faces.
These are obtained by turning all inequalities in the definition of the cone into equalities, except one. For example, one of the one-dimensional faces is given by
$$
R_1 = \{(\alpha_1,\ldots,\alpha_{n+1})\in L(C)\colon \alpha_1+ \sigma \alpha_{n+1} \leq 0, \alpha_2+ \sigma \alpha_{n+1} =0, \ldots,\alpha_n + \sigma \alpha_{n+1} = 0\}.
$$
Taking $\alpha_{n+1} = \sigma/(1+n\sigma^2)$, we obtain that $R_1$ is a ray spanned by the vector
$$
\left(1-\frac{\sigma^2}{1+n\sigma^2}, -\frac{\sigma^2}{1+n\sigma^2}, \ldots, -\frac{\sigma^2}{1+n\sigma^2}, \frac{\sigma}{1+n\sigma^2}\right),
$$
where the value of the first coordinate was computed using the linear relation in the definition of $L(C)$. Thus, the cone $C^{\circ}\cap L(C)$ is spanned by the vectors
\begin{equation}\label{eq:v_i_def}
v_i := e_i -\frac{\sigma^2}{1+n\sigma^2}(e_1+\ldots+ e_n)+ \frac{\sigma}{1+n\sigma^2}e_{n+1} ,
\quad
i=1,\ldots,n.
\end{equation}
It is easy to verify that
\begin{equation}\label{eq:v_scalar_products}
\langle v_i, v_j \rangle =
\begin{cases}
1-\frac{\sigma^2}{1+n\sigma^2} , &\text{ if } i=j,\\
-\frac{\sigma^2}{1+n\sigma^2}, &\text{ if } i\neq j.
\end{cases}
\end{equation}
Thus, the cone spanned by $v_1,\ldots,v_n$ is isometric to $C_{n}(-\frac{\sigma^2}{1+n\sigma^2})$.
\end{proof}

\subsection{Proof of Proposition~\ref{prop:geometric_interpretation_g_n}}
We prove that
$$
\alpha (C_n(r)) = g_n\left(-\frac{r}{1+ n r}\right).
$$
Consider the cone $D\subset \R^N$ spanned by the vectors $v_1,\ldots, v_n$ such that
$$
\langle v_i, v_j \rangle =
\begin{cases}
1-\frac{r}{1+nr} , &\text{ if } i=j,\\
-\frac{r}{1+nr}, &\text{ if } i\neq j.
\end{cases}
$$
Then, $D$ is isometric to $C_n(-\frac{r}{1+ nr})$. The polar cone is given by
$$
D^\circ = \{x\in \R^{N} \colon \langle x, v_1\rangle \leq 0,\ldots, \langle x, v_n\rangle \leq 0\}.
$$
Let $\xi$ be a standard normal random vector on $\R^{N}$. Then, the solid angle of $D^{\circ}$ is given by
$$
\alpha(D^{\circ}) = \P[\xi \in D^{\circ}] =\P\left[ \langle \xi,v_1\rangle \leq 0, \ldots,\langle \xi,v_n\rangle\leq 0 \right].
$$
Introducing the random variables $\eta_i := \langle \xi,v_i\rangle$, $i=1,\ldots,n$,
we observe that the random vector $(\eta_1,\ldots,\eta_n)$ is zero mean Gaussian with covariances given by
\begin{equation}\label{eq:eta_covar1}
\Cov (\eta_i, \eta_j)=
\langle v_i, v_j\rangle
=
\begin{cases}
1-\frac{r}{1+nr} , &\text{ if } i=j,\\
-\frac{r}{1+nr}, &\text{ if } i\neq j.
\end{cases}
\end{equation}
Hence, $\alpha(D^{\circ}) = \P[\eta_1\leq 0,\ldots,\eta_n \leq 0] = g_n\left(-\frac{r}{1+ n r}\right)$ by definition of $g_n$. On the other hand, $D^{\circ}$ is the direct sum of $L(D)\cap D^{\circ}$ and the orthogonal complement of $L(D)$. It follows that
$$
\alpha(L(D)\cap D^{\circ}) = \alpha(D^{\circ}) = g_n\left(-\frac{r}{1+ nr}\right).
$$
By Proposition~\ref{prop:dual_cone}, $L(D)\cap D^{\circ}$ is isometric to $C_n(r)$, thus completing the proof.
\hfill $\Box$

\subsection{Internal and normal angles: Proof of (\ref{eq:angles_claim})}
Recall that $C$ is a cone given by~\eqref{eq:C_def} and that the linear hull of $C$ is a codimension $1$ linear subspace of $\R^n$ given by
\begin{equation}\label{eq:L_C_def}
L(C) =\{(\alpha_1,\ldots,\alpha_n, (\alpha_1+\ldots+\alpha_n) \sigma ) \colon \alpha_1,\ldots,\alpha_n\in\R\}.
\end{equation}
Inside $L(C)$, the convex cone $C$ is defined by the inequalities $\alpha_1\geq 0,\ldots,\alpha_n\geq 0$. The $k$-faces of $C$ are obtained by turning $n-k$ of these inequalities into equalities, therefore the number of $k$-faces is $\binom nk$.  Without restriction of generality, we consider a $k$-face $F$ of the form
\begin{equation}\label{eq:F_def}
F=\{(\alpha_1,\ldots,\alpha_k, \underbrace{0,\ldots,0}_{n-k}, (\alpha_1+\ldots+\alpha_k) \sigma ) \colon \alpha_1,\ldots,\alpha_k \geq 0\}.
\end{equation}
Since $F$ is isometric to $C_k(\sigma^2)$, Proposition~\ref{prop:geometric_interpretation_g_n} yields the following formula for its solid angle:
$$
\alpha (F) = \alpha (C_k(\sigma^2)) = g_k\left(-\frac{\sigma^2}{1+ k \sigma^2}\right).
$$
To compute $\alpha(N_F(C))$, observe that by the polar correspondence, $N_F(C)\cap L(C)$ is some $(n-k)$-dimensional face of the polar cone $C^{\circ}\cap L(C)$. The latter cone is isometric to $C_{n}(-\frac{\sigma^2}{1+n\sigma^2})$ by Proposition~\ref{prop:dual_cone}. Since all $(n-k)$-dimensional faces of $C_{n}(-\frac{\sigma^2}{1+n\sigma^2})$ are isometric to $C_{n-k}(-\frac{\sigma^2}{1+n\sigma^2})$, we can apply Proposition~\ref{prop:geometric_interpretation_g_n} to obtain that
\begin{multline*}
\alpha(N_F(C)) = \alpha(N_F(C)\cap L(C)) = \alpha\left(C_{n-k}\left(-\frac{\sigma^2}{1+n\sigma^2}\right)\right)
\\
= g_{n-k}\left(- \frac{-\frac{\sigma^2}{1+n\sigma^2}}{1-  \frac{(n-k)\sigma^2}{1+n\sigma^2}} \right)
=
g_{n-k} \left(\frac{\sigma^2}{1+k\sigma^2}\right).
\end{multline*}
This completes the proof of~\eqref{eq:angles_claim} and of Theorem~\ref{theo:main1}.
\hfill $\Box$


\subsection{Proof of Proposition~\ref{prop:angles_simplex}}
By symmetry, we may consider the face of the form $F= \conv(e_1,\ldots,e_{k})$. It follows from~\eqref{eq:def_tangent_cone} that the tangent cone is given by
$$
T_F(\Delta_n) = \{(\alpha_1,\ldots,\alpha_n)\in\R^n\colon \alpha_1+\ldots+\alpha_n = 0, \alpha_{k+1}\geq 0,\ldots, \alpha_n\geq 0\}.
$$
Thus, $T_F(\Delta_n)$ is a direct orthogonal sum of the linear subspace $L_{n,k}$ given by $\alpha_1+\ldots+\alpha_k=0$, $\alpha_{k+1}=\ldots=\alpha_n = 0$ (which is the lineality space of the cone) and the cone $D_{n,k} = \pos (u_1,\ldots,u_{n-k})$ spanned by the $n-k$ vectors
$$
u_i := -(e_1+\ldots+e_k)/k + e_{i+k}, \quad i=1,\ldots,n-k.
$$
The scalar products of these vectors are given by
$$
\langle u_i, u_j\rangle
=
\begin{cases}
1+\frac 1k, &\text{ if } i=j,\\
\frac 1k, &\text{ if } i\neq j.
\end{cases}
$$
Hence, the cone $D_{n,k}$ is isometric to $C_{n-k}(1/k)$. From Proposition~\ref{prop:geometric_interpretation_g_n} we deduce that the solid angle of $T_F(\Delta_n)$ is $g_{n-k}(-1/n)$.

The normal cone $N_F(\Delta_n)= T_F^\circ(\Delta_n)$ is the direct orthogonal sum of the line $\{\alpha_1=\ldots=\alpha_n\}$ and the polar cone of $D_{n,k}$ taken w.r.t.\ the ambient space $L(D_{n,k})$. The latter cone is isometric to $C_{n-k}(-1/n)$ by Proposition~\ref{prop:dual_cone}. From Proposition~\ref{prop:geometric_interpretation_g_n} we deduce that the solid angle of $N_F(\Delta_n)$ equals $g_{n-k}(1/k)$.
\hfill $\Box$

\subsection{Proof of Proposition~\ref{prop:upsilon_n_k}}
If $r= \sigma^2\geq 0$, then the proof follows immediately from~\eqref{eq:angles_claim}.
Let $r\in (-\frac 1n, 0)$. For a cone $D\subset \R^N$ we have the relation $\upsilon_k(D) = \upsilon_{N-k}(D^\circ)$; see~\cite[Section~2.2]{amelunxen_lotz}. Applying this relation to the cone $C_n(r)$ in the ambient space $L(C_n(r))$ and recalling Proposition~\ref{prop:dual_cone}, we obtain
$$
\upsilon_k(C_n(r)) = \upsilon_{n-k}\left(C_n\left(-\frac{r}{1+nr}\right)\right)
=
b_{n,n-k}\left(-\frac{r}{1+nr}\right),
$$
where the last step follows from the already established part of Proposition~\ref{prop:upsilon_n_k} and the fact that $-\frac{r}{1+nr}>0$. Using the definition of $b_{n,n-k}$, we obtain
\begin{multline*}
b_{n,n-k}\left(-\frac{r}{1+nr}\right)
=
\binom {n}{n-k}
g_{n-k}\left(- \frac{-\frac{r}{1+nr}}{1-\frac{(n-k)r}{1+nr}}\right) g_k \left(\frac{-\frac{r}{1+nr}}{1 - \frac{(n-k)r}{1+nr}}\right)
\\=
\binom {n}{k}
g_k \left(-\frac{r}{1+kr}\right) g_{n-k}\left(\frac{r}{1+kr}\right)
=
b_{n,k}(r),
\end{multline*}
which proves the claimed formula.
\hfill $\Box$

\section{Properties of \texorpdfstring{$g_n$}{gn}}
\subsection{Proof of Proposition~\ref{prop:g_n_properties}}\label{subsec:g_n_properties}
In the following let $(\eta_1,\ldots,\eta_n)$ be a zero-mean Gaussian vector whose covariance matrix $\Sigma = (r_{ij})_{i,j=1}^n$ given by
$$
r_{ij} = \Cov (\eta_i,\eta_j) =
\begin{cases}
1+r, &\text{ if } i=j,\\
r, &\text{ if } i\neq j.
\end{cases}
$$
Using the inequality between the arithmetic and quadratic means, it is easy to check that this matrix is positive semidefinite for $r\geq -1/n$.
Recall that by definition
\begin{equation}\label{eq:g_n_r}
g_n(r) = \P[\eta_1<0,\ldots,\eta_n <0]= \P[\eta_1>0,\ldots,\eta_n>0], \quad r\geq -1/n.
\end{equation}

\begin{proof}[Proof of (a)]
Let $r>-1/n$ be real. (The case $r=-1/n$ can be then deduced from the continuity of $g_n$ at $-1/n$.)
It is straightforward to check that $\det \Sigma = 1+rn>0$ and that the inverse matrix $\Sigma^{-1} = (s_{ij})_{i,j=1}^n$ is given by
$$
s_{ij} =
\begin{cases}
1 - \frac{r}{1+nr}, &\text{ if } i=j,\\
- \frac{r}{1+nr}, &\text{ if } i\neq j.
\end{cases}
$$
Using~\eqref{eq:g_n_r} and the formula for the multivariate Gaussian density, we obtain for all $r>-1/n$,
\begin{equation}\label{eq:g_n_as_mult_integral}
g_n(r) = \frac{1}{(2\pi)^{n/2} \sqrt{1+nr}}  \int_{(0,\infty)^n}
\exp\left\{ \frac{r}{2(1+nr)} \left(\sum_{j=1}^n x_j\right)^2 - \frac 12\sum_{j=1}^n x_j^2\right\} \dd x_1 \ldots \dd x_n.
\end{equation}
The integral converges for complex $r$ satisfying $\Re \frac{r}{1+nr} <\frac 1n$, which is equivalent to $\Re r > -\frac 1n$. Indeed, by the inequality $(\sum_{j=1}^n x_j)^2 \leq n \sum_{j=1}^n x_j^2$, we have
$$
\Re \left(\frac{r}{2(1+nr)} \left(\sum_{j=1}^n x_j\right)^2 - \frac 12\sum_{j=1}^n x_j^2 \right)
\leq
\frac 12 \left(\Re\frac{r}{(1+nr)} - \frac 1n\right) \left(\sum_{j=1}^n x_j\right)^2.
$$
Hence, the right-hand side of~\eqref{eq:g_n_as_mult_integral} defines an analytic function of $r$ in the half-plane $\Re r > -\frac 1n$.  In particular, $g_n(r)$ has an analytic continuation to this half-plane.

Next we prove that for all $r>0$,
\begin{equation}\label{eq:g_n_integral}
g_n(r) = \frac 1 {\sqrt{2\pi}} \int_0^{\infty} (\Phi^n(\sqrt r x) + \Phi^n(-\sqrt r x)) \eee^{-x^2/2} \dd x.
\end{equation}
Let $\xi,\xi_1,\ldots,\xi_n$ be i.i.d.\ standard Gaussian random variables. We have a distributional representation $\eta_ k= \xi_k - \sqrt r\xi$, $k=1,\ldots,n$. It follows that
\begin{align}
g_n(r)
&=
\P[\eta_1<0,\ldots,\eta_n <0]\notag\\
&=
\P[\xi_1 < \sqrt r \xi, \ldots, \xi_n < \sqrt r \xi]\notag\\
&=
\frac 1 {\sqrt{2\pi}} \int_{-\infty}^{+\infty} \Phi^{n}(\sqrt r x) \eee^{-x^2/2} \dd x,\label{eq:g_n_in_the_proof}
\end{align}
which yields~\eqref{eq:g_n_integral} after splitting the integral.

It remains  to prove that the right-hand side of~\eqref{eq:g_n_integral} is an analytic function of $r$ in the half-plane $\Re r > -\frac 1n$, which would imply that~\eqref{eq:g_n_integral} holds in this half-plane by the uniqueness of the analytic continuation.  First of all, observe that for every fixed $x>0$, the expression $\Phi^n(z x) + \Phi^n(- z x)$ is an analytic function of $z\in \C$ which remains invariant under the substitution  $z\mapsto -z$. Hence, it can be written as an everywhere convergent  Taylor series in even powers of $z$. It follows that for every fixed $x>0$, the expression $(\Phi^n(\sqrt r x) + \Phi^n(-\sqrt r x))$ defines an analytic function of $r\in\C$. To prove the analyticity of the integral on the right-hand side of~\eqref{eq:g_n_integral}, we argue as follows.  We have the estimate
$$
|\Phi(z)| \leq C \max \{1,  |\eee^{-z^2/2}|\} = C \max \{1,  \eee^{-\Re (z^2)/2}\},\quad z\in\C.
$$
It follows that
$$
|\Phi^n(\sqrt r x) + \Phi^n(-\sqrt r x)| \eee^{-x^2/2} \leq 2 C \max \{\eee^{-x^2/2}, \eee^{- \Re (1+ nr) x^2/2}\}.
$$
By the dominated convergence theorem, the right-hand side of~\eqref{eq:g_n_integral} is a continuous function of $r$ on the half-plane $\Re r > -\frac 1n$. Moreover, by Fubini's theorem and by the analyticity of $(\Phi^n(\sqrt r x) + \Phi^n(-\sqrt r x))$, the integral of  the right-hand side of~\eqref{eq:g_n_integral} along any triangular contour  vanishes. By Morera's theorem, the right-hand side of~\eqref{eq:g_n_integral} is an analytic function in the half-plane $\Re r > -\frac 1n$.
By the uniqueness principle for analytic functions, \eqref{eq:g_n_integral} must hold in this half-plane.
\end{proof}

\begin{proof}[Proof of (b)]
By analyticity, it suffices to consider $r>0$.
Differentiating under the sign of the integral in~\eqref{eq:g_n_in_the_proof}, we obtain
\begin{align*}
g_n'(r)
&=
\frac 1{\sqrt{2\pi}} \int_{-\infty}^{+\infty} n \Phi^{n-1}(\sqrt r x) \frac 1 {\sqrt{2\pi}} \eee^{-rx^2/2} \frac x {2\sqrt r} \eee^{-x^2/2} \dd x\\
&=
\frac n {4\pi\sqrt r} \int_{-\infty}^{+\infty}  \Phi^{n-1}(\sqrt r x)   x \eee^{-(r+1)x^2/2} \dd x.
\end{align*}
Writing $x \eee^{-(r+1)x^2/2}\dd x = -\frac 1 {r+1} \dd  \eee^{-(r+1)x^2/2}$ and integrating by parts yields
\begin{align*}
g_n'(r)
&=
\frac n {4\pi\sqrt r} \int_{-\infty}^{+\infty} (n-1) \Phi^{n-2}(\sqrt r x) \frac 1 {\sqrt{2\pi}} \eee^{-rx^2/2} \sqrt r \cdot  \frac 1 {r+1}\eee^{-(r+1)x^2/2} \dd x\\
&=
\frac {n(n-1)} {4\pi (r+1)} \frac 1 {\sqrt{2\pi}} \int_{-\infty}^{+\infty} \Phi^{n-2}(\sqrt r x) \eee^{-(2r+1)x^2/2} \dd x.
\end{align*}
Finally, making the change of variables $y:= \sqrt{2r+1}\, x$, we arrive at
\begin{align*}
g_n'(r)
&= \frac {n(n-1)} {4\pi (r+1)\sqrt{2r+1}} \frac 1 {\sqrt{2\pi}} \int_{-\infty}^{+\infty} \Phi^{n-2}\left(\sqrt {\frac {r}{2r+1}} \, y\right)\eee^{-y^2/2}  \dd y\\
&=
\frac {n(n-1)} {4\pi (r+1)\sqrt{2r+1}} g_{n-2}\left(\frac {r}{2r+1}\right).
\end{align*}
\end{proof}

\begin{proof}[Proof of (c)]
By definition, $g_1(r) = \P[\eta_1 <0] = \frac 12$, since $\eta_1$ is centered Gaussian.
\end{proof}

\begin{proof}[Proof of (d) and (e)]
In the case $r=0$,  the random variables $\eta_1,\ldots,\eta_n$ are independent standard Gaussian and hence
$$
g_n(0) = \P[\eta_1<0,\ldots,\eta_n<0] = 2^{-n}.
$$
In the case $r=1$, we have a distributional representation $\eta_i = \xi_i - \xi$, where $\xi,\xi_1,\ldots,\xi_n$ are independent standard Gaussian. Hence
$$
g_n(1) = \P[\xi_1-\xi <0,\ldots,\xi_n-\xi <0] = \P[\max\{\xi_1,\ldots,\xi_n\} < \xi ] = \frac {1}{n+1}
$$
because any of the values $\xi,\xi_1,\ldots,\xi_n$ can be the maximum with the same probability.

To prove that $\lim_{r\to +\infty} g_n(r) = \frac 12$, use~\eqref{eq:g_n_integral} together with the relation
$$
\lim_{r\to +\infty} (\Phi^n(\sqrt r x) + \Phi^n(-\sqrt r x)) \eee^{-x^2/2} = \eee^{-x^2/2}, \quad x\geq 0,
$$
and the dominated convergence theorem.

Finally, in the case $r=-1/n$ we have the linear relation $\eta_1+\ldots+\eta_n = 0$ (which can be verified by showing that the variance of $\eta_1+\ldots+\eta_n$ vanishes), hence $g_n(-\frac 1n) = 0$.
\end{proof}

\begin{proof}[Proof of (f)]
Let $n=2$. Introduce the variables  $\eta_1^*:= \eta_1/\sqrt{r+1}$ and $\eta_2^* := \eta_2 /\sqrt{r+1}$ which have joint Gaussian law with unit variances and covariance $r/(1+r)$. It follows that
$$
g_2(r) = \P[\eta_1 < 0, \eta_2<0] = \P[\eta_1^* < 0, \eta_2^* <0] = \frac 14 + \frac 1{2\pi} \arcsin \frac{r}{1+r},
$$
by the well-known Sheppard formula for the quadrant probability of a bivariate Gaussian density; see~\cite[p.~121]{bingham_doney} and the references therein.

An alternative proof of this identity is based on Part (b). We give only the proof for $g_3$, since the proof for $g_2$ is similar.  By Part (b), $g_3$ satisfies the differential equation  $g_3'(r) = \frac 3 {4\pi} \frac{1}{(r+1)\sqrt{2r+1}}$   together with the initial condition $g_3(0)= 1/8$. It is easy to check that $g_3(r) = \frac 18 + \frac 3 {4\pi} \arcsin \frac{r}{1+r}$ is the solution.
\end{proof}

\begin{proof}[Proof of (g)]
Using~\eqref{eq:g_n_as_mult_integral} with $r = -\frac 1n + \eps$ and then introducing the variables $y_i := x_i / \sqrt \eps$ yields
\begin{align*}
g_n\left(-\frac 1n + \eps\right)
&= \frac 1 {(2\pi)^{n/2} \sqrt{n\eps}} \int_{(0,\infty)^n} \exp\left\{ \frac {-\frac 1n + \eps}{2\eps n} \left(\sum_{i=1}^n x_i\right)^2 - \frac 12 \sum_{i=1}^n x_i^2 \right\} \dd x_1\ldots \dd x_n\\
&= \frac {\eps^{n/2}} {(2\pi)^{n/2} \sqrt{n\eps}}\int_{(0,\infty)^n} \exp\left\{ \frac {-\frac 1n + \eps}{2n} \left(\sum_{i=1}^n y_i\right)^2 - \frac \eps 2 \sum_{i=1}^n y_i^2 \right\} \dd y_1\ldots \dd y_n\\
&\sim \frac {\eps^{(n-1)/2}} {(2\pi)^{n/2} \sqrt{n}} \int_{(0,\infty)^n} \exp\left\{ -\frac {1}{2n^2} \left(\sum_{i=1}^n y_i\right)^2\right\} \dd y_1\ldots \dd y_n,
\end{align*}
as $\eps \downarrow 0$.
The volume of the simplex $\{y_1,\ldots,y_n\geq 0, y_1+\ldots+ y_n\leq s\}$ is $s^n/n!$. Hence, the integral on the right-hand side equals
$$
\int_0^\infty \eee^{-s^2/(2n^2)} \dd \left(\frac{s^{n}}{n!}\right) = \int_0^\infty \eee^{-s^2/(2n^2)} \frac{s^{n-1}\dd s}{(n-1)!} = 2^{\frac n2-1} n^n \Gamma(n/2)/\Gamma(n),
$$
which completes the proof of~(g).
\end{proof}

\begin{proof}[Proof of (h)]
The functions $g_0(r)=1$ and $g_1(r) = 1/2$ are defined on the whole complex plane. Assume, by induction, that $g_{2n-2}$ and $g_{2n-1}$ are defined as multivalued analytic functions everywhere outside the set $\{-1/k\colon k=1,\ldots, 2n-2\}$. In order to define $g_{2n}$ and $g_{2n+1}$ we use the differential equations from Part (b):
\begin{align*}
g_{2n}'(r)
&= \frac{n(2n-1)}{2\pi (r+1)\sqrt{2r+1}} g_{2n-2} \left(\frac{r}{2r+1}\right),\\
g_{2n+1}'(r)
&=
\frac{n(2n+1)}{2\pi (r+1)\sqrt{2r+1}} g_{2n-1} \left(\frac{r}{2r+1}\right).
\end{align*}
It is easy to check that $\frac{r}{2r+1} \in \{-1/k\colon k=1,\ldots, 2n-2\}$ if and only if $r\in \{-\frac 13, -\frac 14, \ldots, -\frac 1 {2n}\}$.  Hence, the right-hand sides of the differential equations are defined as analytic functions on some unramified cover of $\C\bsl \{-1/k\colon k=1,\ldots, 2n\}$ and we can define $g_{2n}$ and $g_{2n+1}$ by path integration.
\end{proof}

\subsection{Proof of Proposition~\ref{prop:integral_phi}}
For $m\in \boldN_0$ and $n=m+2,m+3,\ldots$ write
$$
I(m,n) := \int_{-\infty}^{+\infty} y^{m} \left(\Phi(iy)\eee^{-y^2/2}\right)^n \dd y.
$$
We have to show that $I(m,n)= 0$. From~\eqref{eq:identity_integral_phi} we already know that $I(0,n) = 0$ for all $n=2,3,\ldots$, which is the basis of our induction. Define also $I(-1,n) = 0$ for $n\in\boldN$. It follows from~\eqref{eq:Phi_def} that
$$
\frac{\dd}{\dd y} \Phi(iy) = \frac{i}{\sqrt{2\pi}} \eee^{y^2/2}.
$$
Performing  partial integration, we can write
\begin{align*}
I(m,n)
&=
-\frac{i \sqrt{2\pi}}{n+1} \int_{-\infty}^{+\infty} y^m \eee^{-(n+1) y^2/2} \dd \Phi^{n+1}(iy)\\
&=
\frac{i \sqrt{2\pi}}{n+1} \int_{-\infty}^{+\infty} \Phi^{n+1}(iy) \left(m y^{m-1}\eee^{-(n+1) y^2/2} - (n+1) y^{m+1} \eee^{-(n+1) y^2/2}\right) \dd y\\
&=
\frac{i \sqrt{2\pi}}{n+1} \left(m I(m-1,n+1) - (n+1)I(m+1,n+1)\right).
\end{align*}
Observe that this identity is true also for $m=0$. Assuming by induction that we proved that $I(l,n) = 0$ for all $l=0,1,\ldots,m$ and $n=l+2,l+3,\ldots$, we obtain from the above identity that $I(m+1,n+1) = 0$.

Let now $n= m+1$, in which case the integral diverges and we have to pass to the Cauchy principal value.  Write
$$
I(m,m+1) := \int_{0}^{+\infty} y^{m} \left(\left(\Phi(iy)\eee^{-y^2/2}\right)^{m+1} + (-1)^m \left(\Phi(-iy)\eee^{-y^2/2}\right)^{m+1}\right) \dd y.
$$
We need to prove that
$$
I(m,m+1) = \sqrt{\frac \pi 2} \left(\frac i {\sqrt{2\pi}}\right)^m.
$$
To treat the case $m=0$, we observe that $\Phi(iy) + \Phi(-iy) = 1$ for $y\in\R$, see~\eqref{eq:Phi_def},  whence
$$
I(0,1) = \int_0^{\infty}  \left(\Phi(iy)\eee^{-y^2/2} + \Phi(-iy)\eee^{-y^2/2}\right) \dd y
=
\int_0^{\infty}  \eee^{-y^2/2} \dd y
=
\sqrt{\frac \pi 2}.
$$
We proceed by induction. Observe that
$$
\frac \dd {\dd y} \left( \Phi^{m+2}(iy) + (-1)^{m+1} \Phi^{m+2}(-iy) \right)
=
(m+2) \left( \Phi^{m+1}(iy) + (-1)^m \Phi^{m+1}(-iy) \right) \frac{i}{\sqrt{2\pi}} \eee^{y^2/2}.
$$
Integration by parts yields
\begin{align*}
\lefteqn{I(m,m+1)}\\
&=
-\frac{i \sqrt{2\pi}}{m+2} \int_{0}^{+\infty} y^m \eee^{-(m+2) y^2/2} \dd \left( \Phi^{m+2}(iy) + (-1)^{m+1} \Phi^{m+2}(-iy) \right) \\
&=
\frac{i \sqrt{2\pi}}{m+2} \int_{0}^{+\infty} \left( \Phi^{m+2}(iy) + (-1)^{m+1} \Phi^{m+2}(-iy) \right) \left(m y^{m-1} - (m+2) y^{m+1} \right) \eee^{-(m+2) y^2/2} \dd y\\
&=
\frac{i \sqrt{2\pi}}{m+2} \left(m I(m-1,m+2) - (m+2)I(m+1,m+2)\right).
\end{align*}
But we already know that $m I(m-1,m+2) = 0$, whence $I(m,m+1) = - i \sqrt{2\pi}I(m+1,m+2)$ for all $m\in\boldN_0$, and the claim follows by induction.
\hfill $\Box$

\section{Inverting the Laplace transforms}

\subsection{Proof of Corollary~\ref{cor:intensity}}
Conditioning on the event $|X| = r$ and noting that $|X|$ has a $\chi$-distribution with $d$ degrees of freedom, we can write
\begin{align*}
\P[\sigma X\notin \conv (X_1,\ldots,X_n)]
&=
\int_{0}^\infty f_{n,d} (\sigma r) \frac{2^{1-(d/2)}}{\Gamma(d/2)} r^{d-1} \eee^{-r^2/2} \dd r\\
&=
\frac{1}{\Gamma(d/2) \sigma^d} \int_{0}^\infty f_{n,d} (\sqrt{2u})   u^{(d/2) -1} \eee^{-u/\sigma^2} \dd u,
\end{align*}
where we made the  change of variables $\sigma r = \sqrt{2 u}$, $\dd r = \sigma^{-1} \dd u / \sqrt{2u}$.  Taking $\lambda:= 1/\sigma^2$ and applying Theorem~\ref{theo:main1}, we obtain
$$
\frac{\lambda^{d/2}}{\Gamma(d/2)} \int_{0}^\infty f_{n,d} (\sqrt{2u})   u^{(d/2) -1} \eee^{- \lambda u} \dd u
=
2(b_{n,d-1} (1/\lambda) + b_{n,d-3} (1/\lambda) + \ldots),
$$
which proves the theorem. \hfill $\Box$

\subsection{Proof of Theorem~\ref{theo:main_intensity}}
From Corollary~\ref{cor:intensity} and Wendel's formula~\eqref{eq:wendel_formula_rep} we know that
\begin{multline*}
\int_0^{\infty} (f_{n,d}(\sqrt{2u}) - f_{n,d}(0)) u^{(d/2) - 1} \eee^{-\lambda u} \dd u
\\=
2\Gamma(d/2) \lambda^{-d/2} \left(b_{n,d-1}(1/\lambda) - \frac 1 {2^n}\binom {n}{d-1} + b_{n,d-3}(1/\lambda) - \frac 1 {2^n}\binom{n}{d-3} +\ldots\right).
\end{multline*}
By the uniqueness of the Laplace transform, it suffices to show that
$$
\int_0^{\infty} a_{n,k}(u)  \eee^{-\lambda u} \dd u =  \Gamma(d/2) \lambda^{-d/2} \left(b_{n,k}(1/\lambda) - \frac 1 {2^n}\binom {n}{k}\right).
$$
Recalling the formulas for $a_{n,k}(u)$ and $b_{n,k}(1/\lambda)$, see~\eqref{eq:def_b_k}, \eqref{eq:a_k}, we rewrite this as
\begin{multline*}
\int_0^{\infty} \eee^{-\lambda u} \frac 1 {\Gamma(d/2)} \left(\int_0^u \eee^{-vk} F'_{k,n-k} (v) (u-v)^{(d/2)-1}  \dd v\right) \dd u
\\
= \lambda^{-d/2} \left( g_k\left(-\frac{1}{\lambda+k}\right) g_{n-k} \left(\frac{1}{\lambda+k}\right)- \frac 1 {2^n}\right).
\end{multline*}
The inner integral on the left-hand side is the fractional Riemann--Liouville integral of order $d/2$ of the function $\eee^{-vk} F'_{k,n-k} (v)$. Recall that fractional integral of order $\alpha\geq 0$ is defined by
$$
J_\alpha f (u)  = \frac 1 {\Gamma(\alpha)} \int_0^u f(v) (u-v)^{\alpha-1} \dd v
$$
and its Laplace transform is just $\lambda^{-\alpha}$ times the Laplace transform of $f$:
$$
\int_{0}^{\infty} J_\alpha f (u) \eee^{-\lambda u} \dd u = \lambda^{-\alpha}\int_{0}^{\infty} f (u) \eee^{-\lambda u} \dd u.
$$
Using this property, we deduce that it is sufficient to prove that
$$
\int_0^{\infty} \eee^{-\lambda u}  \eee^{-uk} F'_{k,n-k} (u)  \dd u = g_k\left(-\frac{1}{\lambda+k}\right) g_{n-k} \left(\frac{1}{\lambda+k}\right)- \frac 1 {2^n}.
$$
Writing $\mu:= \lambda + k$, we rewrite this as
$$
\int_0^{\infty} \eee^{-\mu u}  F'_{k,n-k} (u)  \dd u = g_k(-1/{\mu}) g_{n-k} (1/\mu) - \frac 1 {2^n}.
$$
Observe that $\lim_{u\downarrow 0}F_{k,n-k}(u) = 2^{-n}$. This follows from~\eqref{eq:F_k_n_minus_k} by observing that
$$
\Phi^{n-k}(0) + \Phi^{n-k}(0) = \frac 2 {2^{n-k}},
\quad
\Re \Phi^{k} (0) = \frac 1 {2^k},
\quad
\int_{0}^1 \frac{\dd w}{\sqrt{w(1-w)}} 
=\pi
$$
and using dominated convergence. Using partial integration, we can write the above as
\begin{equation}\label{eq:need1}
\mu \int_0^{\infty} \eee^{-\mu u} F_{k,n-k} (u)  \dd u = g_k(-1/{\mu}) g_{n-k} (1/\mu).
\end{equation}
Recall from~\eqref{eq:F_k_n_minus_k} that $F_{k,n-k}(u) = \int_{0}^{u} h_{n-k}^{(1)} (w) h_k^{(2)} (u-w) \dd w$, where
$$
h_{n-k}^{(1)}(w) = \frac{\Phi^{n-k}(\sqrt{2 w}) + \Phi^{n-k}(-\sqrt{2 w})}{2\sqrt {\pi w}},
\quad
h_{k}^{(2)}(w) = \frac{\Phi^{k}(i \sqrt{2 w}) + \Phi^{k}(-i\sqrt{2 w})}{2\sqrt {\pi w}}.
$$
Let us compute the Laplace transforms of $h_{n-k}^{(1)}$ and $h_{k}^{(2)}$:
\begin{align*}
\int_0^{\infty} \eee^{-\mu w} h_{n-k}^{(1)}(w) \dd w
&=
\int_0^{\infty} \eee^{-\mu w} \left(\Phi^{n-k}(\sqrt{2 w}) + \Phi^{n-k}(-\sqrt{2 w})\right) \frac{\dd w}{2\sqrt {\pi w}}\\
&=
\frac 1 {\sqrt{2\pi \mu}} \int_0^{\infty} \eee^{-y^2/2} \left(\Phi^{n-k}\left(\frac{y}{\sqrt \mu}\right) + \Phi^{n-k}\left(-\frac{y}{\sqrt \mu}\right)\right)  \dd y\\
&=
\mu^{-1/2} g_{n-k}(1/\mu),
\end{align*}
where we made the change of variables $w = y^2/(2\mu)$ in the second step and recalled~\eqref{eq:g_n_r_part_a} in the last step.
Arguing in an analogous way, we obtain
\begin{align*}
\int_0^{\infty} \eee^{-\mu w} h_{k}^{(2)}(w) \dd w
&=
\int_0^{\infty} \eee^{-\mu w} \left(\Phi^{k}(i\sqrt{2 w}) + \Phi^{k}(-i\sqrt{2 w})\right) \frac{\dd w}{2\sqrt {\pi w}}\\
&=
\frac 1 {\sqrt{2\pi \mu}} \int_0^{\infty} \eee^{-y^2/2} \left(\Phi^{k}\left(i\frac{y}{\sqrt \mu}\right) + \Phi^{k}\left(-i\frac{y}{\sqrt \mu}\right)\right)  \dd y\\
&=
\mu^{-1/2} g_{k}(-1/\mu).
\end{align*}
Since the Laplace transform of the convolution is the product of the Laplace transforms, we arrive at~\eqref{eq:need1}, which completes the proof.
\hfill $\Box$

\subsection{Proof of Theorem~\ref{theo:intensity_d_2}}
Corollary~\ref{cor:intensity} with $d=2$ states that for all $\lambda>0$,
$$
\int_0^{\infty} f_{n,2}(\sqrt{2u}) \eee^{-\lambda u} \dd u
= \frac 2 {\lambda} b_{n,1}(1/\lambda) = \frac {n}{\lambda} g_{n-1}\left(\frac{1}{1+\lambda}\right),
$$
where we used~\eqref{eq:def_b_k} and the fact that $g_1(r) = 1/2$. Recalling the formula for $g_{n-1}$, see~\eqref{eq:g_n_r_part_a}, we arrive at
\begin{align*}
\int_0^{\infty} f_{n,2}(\sqrt{2u}) \eee^{-\lambda u} \dd u
&=
\frac {n}{\lambda\sqrt{2\pi}} \int_0^{\infty} \left(\Phi^{n-1}\left(\frac {x}{\sqrt{1+\lambda}}\right) + \Phi^{n-1}\left(-\frac {x}{\sqrt{1+\lambda}}\right)\right) \eee^{-x^2/2}\dd x\\
&=
\frac{\sqrt{1+\lambda}}{\lambda} \cdot \frac {n}{\sqrt{2\pi}} \int_0^{\infty} \frac{\Phi^{n-1}(\sqrt{2z}) + \Phi^{n-1}(-\sqrt{2z})}{\sqrt{2z}} \eee^{-z} \eee^{-z\lambda}\dd z\\
&=
\frac{\sqrt{1+\lambda}}{\lambda} \cdot \int_0^{\infty} f_{M_n^2/2}(z) \eee^{-z\lambda} \dd z,
\end{align*}
where we used the change of variables $x/\sqrt {1+\lambda} = \sqrt{2z}$ and $f_{M_n^2/2}$ denotes the density of the random variable $\frac 12 M_n^2$, namely
$$
f_{M_n^2/2}(z) = \frac{\dd}{\dd z} \left(\Phi^{n}(\sqrt {2z}) - \Phi^{n}(-\sqrt {2z})\right) =
\frac {n}{\sqrt{2\pi}} \frac{\Phi^{n-1}(\sqrt{2z}) + \Phi^{n-1}(-\sqrt{2z})}{\sqrt{2z}} \eee^{-z},
$$
where $z>0$.
The inverse Laplace transform of $\sqrt{1+\lambda} /\lambda$ is given by $(2\Phi(\sqrt{2z}) - 1) + \eee^{-z}/\sqrt{\pi z}$. Observe that the first summand is just the distribution function $F_{\xi^2/2}$ of $\frac 12 \xi^2$ (with $\xi$ standard normal), whereas the second summand is the density $f_{\xi^2/2} = F_{\xi^2/2}'$ of the same random variable. Since the inverse Laplace transform of a product is a convolution of the inverse Laplace transforms, we arrive at
\begin{align*}
f_{n,2}(\sqrt{2u})
&= \int_{0}^{u} f_{M_n^2/2} (z) F_{\xi^2/2}(u-z) \dd z +  \int_{0}^{u} f_{M_n^2/2} (z) f_{\xi^2/2}(u-z) \dd z\\
&=
\P[ M_{n}^2 + \xi^2 \leq 2u] + \frac {\dd}{\dd u}\P[ M_{n}^2 + \xi^2 \leq 2u],
\end{align*}
which proves the first formula in Theorem~\ref{theo:intensity_d_2}. To verify the second formula, observe that by the change of variables $x:= \sqrt {2z}$ and then $y: = x/\sqrt {2u}$,
\begin{align*}
\int_{0}^{u} f_{M_n^2/2} (z) f_{\xi^2/2}(u-z) \dd z
&=
\frac {n}{\sqrt{2\pi}} \int_0^u  \frac{\Phi^{n-1}(\sqrt{2z}) + \Phi^{n-1}(-\sqrt{2z})}{\sqrt{2z}} \eee^{-z} \frac{\eee^{-(u-z)}}{\sqrt \pi \sqrt{u-z}} \dd z\\
&=
n\eee^{-u} \int_{0}^{\sqrt{2u}}  \frac{\Phi^{n-1}(x) + \Phi^{n-1}(-x)}{\pi \sqrt{2u-x^2}} \dd x\\
&=
n\eee^{-u} \int_{-\sqrt{2u}}^{\sqrt{2u}} \frac{\Phi^{n-1}(x) }{\pi \sqrt{2u-x^2}} \dd x\\
&=
n\eee^{-u} \int_{-1}^{1} \frac{\Phi^{n-1}(\sqrt {2u} y) }{\pi \sqrt{1-y^2}} \dd y.
\end{align*}
The integral equals $\P[M_{n-1}\leq \sqrt{2u} W]$ since $\Phi^{n-1}$ is the distribution function of $M_{n-1}$.
\hfill $\Box$

\section{Proof of Theorem~\ref{theo:goodman_pollack}}\label{sec:proof_goodman_pollack}
The proof relies strongly on the ideas of \citet{baryshnikov_vitale} combined with the Bartlett decomposition of the Gaussian matrix and Theorem~\ref{theo:main1}.

\vspace*{2mm}
\noindent
\textit{Step 1: Bartlett decomposition of a Gaussian matrix.}
Let $O$ be a random orthogonal matrix distributed according to the Haar probability measure on the group $O(n)$. Independently of $O$, let $L$ be a random lower triangular $n\times n$-matrix with a.s.\ positive entries on the diagonal and the following distribution. The entries of $L$ are independent, the squared diagonal entries have $\chi^2$-distributions with $n,n-1,\ldots,1$ degrees of freedom, whereas the entries below the diagonal are standard normal.
Define an $n\times n$-matrix $G$ by
$$
G = LO.
$$
It is known (Bartlett decomposition) that the entries of $G$ are independent  standard Gaussian random variables.

\vspace*{2mm}
\noindent
\textit{Step 2: Relating the Goodman--Pollack model to the Gaussian polytope.}
Consider an $n\times (n+1)$-matrix $S$ whose columns are the vectors $v_1,\ldots,v_{n+1}$. Recall also that $\Pi$ is the $d\times n$-matrix of the projection from $\R^n$ onto the linear subspace $\R^d$ spanned by $e_1,\ldots,e_d$. Note that $\Pi$ consists of an identity matrix $I_d$ extended by a zero $d\times (n-d)$ matrix.  Following \citet{baryshnikov_vitale}, consider the $d\times (n+1)$-matrix
\begin{equation}\label{eq:tilde_S}
\tilde S := \Pi G S =  \Pi L O S = \tilde L \Pi O S,
\end{equation}
where $\tilde L$ is a lower-triangular $d\times d$-matrix obtained from $L$ by removing all rows and columns except the first $d$ ones.
The last equality follows from the simple identity $\Pi L = \tilde L \Pi$. It follows from the corresponding properties of $L$ that the matrix $\tilde L$ is lower-triangular, its entries are independent, the squared diagonal entries have $\chi^2$-distributions  with $n,n-1,\ldots, n-d$ degrees of freedom, while the entries below the diagonal are standard Gaussian.

Observe that the matrices $\tilde L$ and $\Pi O S$ are stochastically independent and the columns of the latter matrix are the vectors $\Pi O v_1, \ldots, \Pi O v_{n+1}$ whose convex hull is the Goodman--Pollack polytope $\cQ_{n+1,d}$.

Denote the columns of the matrix $\tilde S$ by $Y_{\bullet, 1},\ldots, Y_{\bullet, (n+1)}$. Let us write $Y_{\bullet, j} = (Y_{1,j}, \ldots, Y_{d,j})^\top$, $1\leq j\leq n+1$. It follows from the formula $\tilde S = \Pi G S$ that the random variables $Y_{i,j}$, $1\leq i\leq d$, $1\leq j\leq n+1$, are jointly Gaussian with mean zero. Let us compute their covariances.  Writing  $G = (g_{i,j})_{1\leq i \leq n, 1\leq j\leq n}$ with independent standard Gaussian entries $g_{i,j}$, we observe that $\Pi G = (g_{i,j})_{1\leq i \leq d, 1\leq j \leq n}$ and hence,
\begin{align*}
\boldE [Y_{i_1, j_1} Y_{i_2, j_2}]
&=
\boldE \left[\langle (g_{i_1, 1}, \ldots, g_{i_1, n})^\top, v_{j_1} \rangle  \langle (g_{i_2, 1}, \ldots, g_{i_2, n})^\top, v_{j_2} \rangle  \right]\\
&=
\begin{cases}
0, &\text{ if }  i_1\neq i_2,\\
\langle v_{j_1}, v_{j_2} \rangle , &\text{ if }  i_1 = i_2.
\end{cases}
\end{align*}
By the properties of the regular simplex we have $\langle v_{j_1}, v_{j_2}\rangle = -1/n <0$ provided that $j_1\neq j_2$.
Let $W$ be a standard Gaussian random vector on $\R^d$ independent of $L,O$ (and hence, $G$ and $\tilde S$).  It follows that
$$
\sqrt{\frac{n}{n+1}} Y_{\bullet, 1} + \frac W{\sqrt{n+1}} , \ldots, \sqrt{\frac{n}{n+1}} Y_{\bullet, n+1} +  \frac W{\sqrt{n+1}}
$$
are mutually independent standard Gaussian random vectors on $\R^d$. Their convex hull has the same distribution as the Gaussian polytope and will be denoted by $\cP_{n+1,d}$ for this reason. Summarizing everything, we obtain the  identity
\begin{multline}\label{eq:baryshnikov_vitale}
\tilde L \cQ_{n+1,d}
=
 \conv(\tilde L \Pi O v_1,\ldots, \tilde L \Pi O v_{n+1})\\
=
 \conv(Y_{\bullet, 1},\ldots, Y_{\bullet, (n+1)}) = \sqrt{\frac{n+1}{n}} \cP_{n+1,d} - \frac W{ \sqrt{n}},
\end{multline}
where $\tilde L$, $W$, $\cQ_{n+1,d}$ are independent. The fact that $\cP_{n+1,d}$ can be obtained from $\cQ_{n+1,d}$ by a random affine transformation stochastically independent from $\cQ_{n+1,d}$ was proved by~\citet{baryshnikov_vitale}. We rederived it because we shall need the explicit form of the affine transformation in what follows.

\vspace*{2mm}
\noindent
\textit{Step 3: Relating  the non-absorption probabilities for $\cP_{n+1,d}$ and $\cQ_{n+1,d}$.}
Let now $X$ be a standard Gaussian vector on $\R^d$ which is independent of everything else. We know from Theorem~\ref{theo:main1} that for every $\sigma^2\geq 0$,
$$
\P[\sigma X \notin \cP_{n+1,d}] = 2(b_{n+1,d-1}(\sigma^2) + b_{n+1,d-3}(\sigma^2) +\ldots).
$$
It follows from~\eqref{eq:baryshnikov_vitale} that
$$
\P\left[ \sqrt{\frac{n+1}{n}} \sigma X  - \frac W{ \sqrt{n}} \notin  \tilde L \cQ_{n+1,d}\right] =  2(b_{n+1,d-1}(\sigma^2) + b_{n+1,d-3}(\sigma^2) +\ldots).
$$
Introducing the standard normal $d$-dimensional vector
$$
\xi :=  \left(\frac{n+1}{n}\sigma^2 + \frac 1n\right)^{-1/2} \left(\sqrt{\frac{n+1}{n}} \sigma X  - \frac W{ \sqrt{n}}\right),
$$
(which is independent of $\tilde L$ and $\cQ_{n+1,d}$) we can rewrite the above as
$$
\P\left[ \sqrt{\frac{n+1}{n}\sigma^2 + \frac 1n} \cdot  \xi  \notin  \tilde L \cQ_{n+1,d}\right] =  2(b_{n+1,d-1}(\sigma^2) + b_{n+1,d-3}(\sigma^2) +\ldots).
$$

Let $Q\in O(d)$ be a deterministic orthogonal matrix.  We claim that the random set $\tilde L \cQ_{n+1,d}$ is invariant w.r.t.\ orthogonal transformations, namely
\begin{equation}\label{eq:tilde_L_Q_invar}
Q \tilde L \cQ_{n+1,d} \eqdistr \tilde L \cQ_{n+1,d}.
\end{equation}
Since $\cQ_{n+1,d}$ is defined as the convex hull of the columns of the matrix $\Pi O S$, it suffices to show that $Q \tilde L \Pi O S \eqdistr \tilde L \Pi O S$, or, equivalently,  $Q \Pi G S \eqdistr \Pi G S$; see~\eqref{eq:tilde_S}. Let $Q'\in O(n)$ be a natural  extension of $Q$ from $\R^d$ to $\R^n$ defined by
$$
Q' e_1= Q e_1,\ldots, Q'e_d = Q e_d, \quad  Q' e_{d+1}=e_{d+1}, \ldots, Q' e_{n}=e_n.
$$
Then, $Q \Pi = \Pi Q'$ and hence, it suffices to show that $\Pi Q' G S \eqdistr \Pi G S$. However, since the entries of $G$ are independent standard Gaussian and the matrix $Q'$ is orthogonal, it is easy to check that $Q' G \eqdistr G$, thus proving~\eqref{eq:tilde_L_Q_invar}.

Let $R_1 := |\xi|$, so that $R^2_1$ has $\chi^2$-distribution with $d$ degrees of freedom. By the orthogonal invariance of the random set $\tilde L \cQ_{n+1,d}$, we can replace $\xi$ by, say, $R_1 e_{n-d}$, thus obtaining
$$
\P\left[ \sqrt{\frac{n+1}{n}\sigma^2 + \frac 1n} \cdot R_1 e_{n-d}  \notin  \tilde L \cQ_{n+1,d}\right] =  2(b_{n+1,d-1}(\sigma^2) + b_{n+1,d-3}(\sigma^2) +\ldots).
$$
Observe that $R_2 := |\tilde L^{-1} e_{n-d}|^{-1}$ satisfies $R_2^2 \sim \chi_{n-d+1}^2$ because of the structure of the lower triangular matrix $\tilde L$.  The random set $\cQ_{n+1,d}$ is also orthogonally invariant (which follows from its definition), hence
$$
\P\left[\sqrt{\frac{n+1}{n}\sigma^2 + \frac 1n} \cdot \frac{R_1}{R_2} e_1  \notin  \cQ_{n+1,d}\right] =  2(b_{n+1,d-1}(\sigma^2) + b_{n+1,d-3}(\sigma^2) +\ldots).
$$
The random variables $R_1\sim \chi_d$ and $R_2\sim \chi_{n-d+1}$ are independent (because $\xi$ and $\tilde L$ are independent), hence the density of $R_1/R_2$ is given by
$$
h(u) = \frac{2 u^{d-1} (1+u^2)^{-(n+1)/2}}{B\left(\frac d2, \frac{n-d+1}{2}\right)}, \quad u\geq 0.
$$
We can finally rewrite the above as
$$
\int_{0}^{\infty} h(u) f_{n+1,d}^* \left(\sqrt{\frac{n+1}{n}\sigma^2 + \frac 1n} \cdot u\right) \dd u = 2(b_{n+1,d-1}(\sigma^2) + b_{n+1,d-3}(\sigma^2) +\ldots),
$$
and Theorem~\ref{theo:goodman_pollack} follows after replacing $n+1$ by $n$.
\hfill $\Box$

\begin{rem}\label{rem:goodman_pollack_proof}
If instead of the regular simplex with $n+1$ vertices inscribed into $\mathbb S^{n-1}$ we rotate and project the regular simplex $\conv (e_1,\ldots,e_n)$, the proof simplifies. The random vectors $Y_{\bullet, 1}, \ldots, Y_{\bullet, n}$ (their number is $n$ rather than $n+1$) are independent standard Gaussian and there is no need of introducing $W$.
\end{rem}

\subsection*{Acknowledgement} Z.\ K.\ is grateful to Alexander Marynych for suggesting to use partial integration in the proof of Proposition~\ref{prop:integral_phi}.

\bibliography{gauss_polytope_bib}
\bibliographystyle{plainnat}

\end{document}